\documentclass[12pt,a4paper]{article}
\usepackage{latexsym,epsfig,fullpage,amsfonts,amssymb,amsmath,amscd,epic,graphics,rotating,theorem,mathabx}
\theoremstyle{plain}
\newtheorem{lemma}{Lemma}

\newtheorem{remark}[lemma]{Remark}

\newtheorem{theorem}{Theorem}

\newtheorem{definition}[lemma]{Definition}

{\theorembodyfont{\rmfamily} 
\font\ncsc=cmcsc10  \font\ntt=cmtt12
\begin{document}
\baselineskip=15pt
\newcommand{\pperp}{\hbox{$\perp\hskip-6pt\perp$}}
\newcommand{\ssim}{\hbox{$\hskip-2pt\sim$}}
\newcommand{\N}{{\mathbb N}}\newcommand{\Delp}{{\Pi}}
\newcommand{\A}{{\mathbb A}}
\newcommand{\Z}{{\mathbb Z}}
\newcommand{\R}{{\mathbb R}}
\newcommand{\C}{{\mathbb C}}\newcommand{\fm}{{\mathfrak{m}}}
\newcommand{\oGam}{{\overline\Gam}}\newcommand{\oh}{{\overline h}}
\newcommand{\Q}{{\mathbb Q}}\newcommand{\T}{{\mathbb T}}\newcommand{\K}{{\mathbb K}}
\newcommand{\PP}{{\mathbb P}}
\newcommand{\mt}{{\operatorname{mt}}}\newcommand{\Cl}{{\operatorname{Cl}}}
\newcommand{\mr}{{\operatorname{m}}}\newcommand{\Trop}{{\operatorname{Trop}}}
\newcommand{\dmp}{{\operatorname{dm}}}
\newcommand{\tp}{{\operatorname{t}}}\newcommand{\ini}{{\operatorname{ini}}}
\newcommand{\dtp}{{\operatorname{dt}}}
\newcommand{\mnote}{\marginpar}\newcommand{\Ev}{{\operatorname{Ev}}}
\newcommand{\Id}{{\operatorname{Id}}}\newcommand{\irr}{{\operatorname{irr}}}\newcommand{\nod}{{\operatorname{nod}}}
\newcommand{\oeps}{{\overline\eps}}\newcommand{\Area}{{\operatorname{Area}}}\newcommand{\End}{{\operatorname{End}}}
\newcommand{\oDel}{{\widetilde\Del}}
\newcommand{\real}{{\operatorname{Re}}}
\newcommand{\conv}{{\operatorname{conv}}}
\newcommand{\Span}{{\operatorname{Span}}}
\newcommand{\Ker}{{\operatorname{Ker}}}
\newcommand{\Fix}{{\operatorname{Fix}}}
\newcommand{\sign}{{\operatorname{sign}}}
\newcommand{\Log}{{\operatorname{Log}}}\newcommand{\dist}{{\operatorname{dist}}}
\newcommand{\oi}{{\overline i}}
\newcommand{\oj}{{\overline j}}
\newcommand{\ob}{{\overline b}}
\newcommand{\os}{{\overline s}}
\newcommand{\oa}{{\overline a}}
\newcommand{\oy}{{\overline y}}
\newcommand{\ow}{{\overline w}}
\newcommand{\ou}{{\overline u}}
\newcommand{\ot}{{\overline t}}
\newcommand{\oz}{{\overline z}}\newcommand{\bu}{{\boldsymbol u}}
\newcommand{\bw}{{\boldsymbol w}}\newcommand{\bi}{{\boldsymbol i}}
\newcommand{\bx}{{\boldsymbol x}}\newcommand{\bp}{{\boldsymbol p}}
\newcommand{\bpp}{{\boldsymbol Q}}\newcommand{\bq}{{\boldsymbol q}}
\newcommand{\by}{{\boldsymbol y}}\newcommand{\bxc}{{\boldsymbol X}}\newcommand{\ba}{{\boldsymbol a}}
\newcommand{\bz}{{\boldsymbol z}}\newcommand{\bv}{{\boldsymbol v}}
\newcommand{\eps}{{\varepsilon}}
\newcommand{\proofend}{\hfill$\Box$\bigskip}
\newcommand{\Int}{{\operatorname{Int}}}
\newcommand{\pr}{{\operatorname{pr}}}
\newcommand{\grad}{{\operatorname{grad}}}
\newcommand{\rk}{{\operatorname{rk}}}
\newcommand{\im}{{\operatorname{Im}}}
\newcommand{\sk}{{\operatorname{sk}}}
\newcommand{\const}{{\operatorname{const}}}
\newcommand{\Sing}{{\operatorname{Sing}}}
\newcommand{\conj}{c}
\newcommand{\Pic}{{\operatorname{Pic}}}
\newcommand{\Crit}{{\operatorname{Crit}}}
\newcommand{\Ch}{{\operatorname{Ch}}}
\newcommand{\discr}{{\operatorname{discr}}}
\newcommand{\Tor}{{\operatorname{Tor}}}
\newcommand{\Conj}{{\operatorname{Conj}}}
\newcommand{\val}{{\operatorname{Val}}}
\newcommand{\res}{{\operatorname{res}}}
\newcommand{\add}{{\operatorname{add}}}
\newcommand{\tmu}{{\C\mu}}
\newcommand{\ov}{{\overline v}}\newcommand{\on}{{\overline n}}
\newcommand{\ox}{{\overline{x}}}
\newcommand{\tet}{{\theta}}
\newcommand{\Del}{{\Delta}}
\newcommand{\bet}{{\beta}}
\newcommand{\kap}{{\kappa}}
\newcommand{\del}{{\delta}}
\newcommand{\sig}{{\sigma}}
\newcommand{\alp}{{\alpha}}
\newcommand{\Sig}{{\Sigma}}
\newcommand{\Gam}{{\Gamma}}
\newcommand{\gam}{{\gamma}}
\newcommand{\Lam}{{\Lambda}}
\newcommand{\lam}{{\lambda}}
\newcommand{\SC}{{SC}}
\newcommand{\MC}{{MC}}
\newcommand{\nek}{{,...,}}
\newcommand{\cim}{{c_{\mbox{\rm im}}}}
\newcommand{\mathto}{\mathop{\to}}
\newcommand{\op}{{\overline p}}

\newcommand{\w}{{\omega}}

\title{Tropical and algebraic curves with multiple points}
\author{Eugenii Shustin}
\date{}
\maketitle

\centerline{\it\small To Oleg Yanovich Viro in occasion of his 60th
birthday}

\begin{abstract}
Patchworking theorems serve as a basic element of the
correspondence between tropical and algebraic curves, which is a
core of the tropical enumerative geometry. We present a new
version of a patchworking theorem which relates plane tropical
curves with complex and real algebraic curves having prescribed
multiple points. It can be used to compute Welschinger invariants
of non-toric Del Pezzo surfaces.
\end{abstract}

\section{Introduction}\label{intro}

The patchworking construction in the toric context is originated
in the Viro method suggested in 1979-80 for obtaining real
algebraic hypersurfaces with prescribed topology
\cite{Vi2,Vi3,Vi4}. Later it was developed and applied to other
problems, in particular, to the tropical geometry. Namely, it
serves as in important step in the proof of a correspondence
between tropical and algebraic curves which in turn is a core of
the enumerative applications of the tropical geometry (see, for
example, the foundational Mikhalkin's work \cite{Mi} and other
versions and modifications in \cite{NS,Sh0,Sh1,Sh2}). We continue
the latter line and present here a new patchworking theorem. The
novelty of our version is that it allows one to patchwork
algebraic curves with prescribed multiple points, whereas the
similar existing statements in tropical geometry apply only to
nonsingular or nodal curves.\footnote{Rephrasing Selman Akbulut,
who called Viro's disciples ``little Viro's", our contribution is
``a little patchworking theorem" descending from ``the great
Viro's patchworking theorem".}

The cited results are restricted to the case of curves in toric
varieties (for example, the plane blown up in at most three
points). Since the consideration of curves on a blown up surface
is equivalent to the study of curves with fixed multiple points on
the original surface, one can apply the tropical enumerative
geometry to count curves on the plane blown up at more than three
points. This approach naturally leads to the question: {\it What
are the plane tropical curves which correspond (as non-Archimedean
amoebas or logarithmic limits) to algebraic curves with fixed
generic multiple points on toric surfaces?} The question appears
to be more complicated than that resolved in \cite{Mi,Sh0}, and no
general answer is known so far.

The goal of the present paper is to prove a patchworking theorem for
a specific sort of plane tropical curves, {\it i.e.} we show that
each tropical curve in the chosen class gives rise to an explicitly
described set of algebraic curves on a given toric surface, in a
given linear system, of a given genus, and with a given collection
of fixed points with prescribed multiplicities (Theorem \ref{t1},
section \ref{sec3}). Furthermore, in the real situation, we compute
the contribution of the constructed curves to the Welschinger
invariant (Theorem \ref{t2}, section \ref{sec3}).

In fact, we do not know all the tropical curves, which may give rise
to the above algebraic curves, and, furthermore, we restrict our
patchworking theorem to a statement which is sufficient to settle
the two following problems:
\begin{itemize}\item to prove recursive formulas of the
Caporaso-Harris type for the Welschinger invariants of
$(\PP^1)^2_{(0,2)}$, the quadric hyperboloid, blown up at two
imaginary points, and for $\PP^2_{(k,2l)}$, $k+2l\le 5$, $l\le 1$,
the plane, blown up at $k$ generic real points and at $l$ pairs of
conjugate imaginary points \cite{IKS};
\item to establish a new correspondence theorem between algebraic
curves of a given genus in a given linear system on a toric surface
and some tropical curves, and find new real tropical enumerative
invariants of real toric surfaces \cite{Sh08}.
\end{itemize} We mention here an important consequence of the
former result

\begin{theorem} {\rm (\cite{IKS})}
Let $\Sig$ be one of the real Del Pezzo surfaces $(\PP^1)^2_{(0,2)}$
or $\PP^2_{(k,2l)}$, $k+2l\le 5$, $l\le 1$, and let $D\subset\Sig$
be a real ample divisor. Then the Welschinger invariants
$W_0(\Sig,D)$, corresponding to the totally real configurations of
points, are positive, and they satisfy the asymptotic relation
$$\lim_{n\to\infty}\frac{\log W_0(\Sig,nD)}{n\log
n}=\lim_{n\to\infty}\frac{\log GW_0(\Sig,nD)}{n\log n}=-K_\Sig D\
,$$ where $GW_0(\Sig,D)$ are the genus zero Gromov-Witten
invariants.
\end{theorem}

A similar statement for all the real toric Del Pezzo surfaces except
for $(\PP^1)^2_{(0,2)}$ was known before \cite{IKS2,IKS3}.

\bigskip

{\bf Preliminary notations and definitions.} If $P\subset\R^n$ is a
pure-dimensional lattice polyhedral complex, $\dim P=d\le n$, by
$|P|$ we denote the lattice volume of $P$, counted so that the
lattice volume of a $d$-dimensional lattice polytope
$\Del\subset\R^n$ is the ratio of the Euclidean volume of $\Del$ and
of the minimal Euclidean volume of a $d$-dimensional lattice simplex
in the linear $d$-subspace of $\R^n$, parallel to the affine
$d$-space spanned by $\Del$. In particular, $|P|=\#P$ if $P$ is
finite.

Given a lattice polyhedron $\Del$, by $\Tor_K(\Del)$\ \footnote{We
omit subindex in the complex case writing simply $\Tor(\Del)$.} we
denote the toric variety over a field $K$, associated with $\Del$,
and by ${\cal L}_\Del$ we denote the tautological line bundle ({\it
i.e.} the bundle generated by the monomials $z^\omega$,
$\omega\in\Del$, as global sections). The divisors
$\Tor_K(\sig)\subset\Tor_K(\Del)$, corresponding to the facets
(faces of codimension $1$) $\sig$ of $\Del$, we call {\bf toric
divisors}. By $\Tor_K(\partial\Del)$ we denote the union of all the
toric divisors in $\Tor_K(\Del)$.

The main field we use is $\K=\bigcup_{m\ge 1}\C((t^{1/m}))$, the
field of locally convergent complex Puiseux series possessing a
non-Archimedean valuation
$$\val:\K^*\to\R,\quad\val\left(\sum_ra_rt^r\right)=-\min\{r\ :\ a_r\ne 0\}\
.$$ Denote
$$\ini\left(\sum_ra_rt^r\right)=a_v,\quad\text{where}\quad
v=-\val\left(\sum_ra_rt^r\right)\ .$$ The field $\K$ is
algebraically closed and contains a closed real subfield
$\K_\R=\Fix(\Conj)$, $\Conj(\sum_ra_rt^r)=\sum_r\overline a_rt^r$.

We remind here the definition of Welschinger invariants \cite{W1},
restricting ourselves to a particular situation. Let $\Sig$ be a
real {\bf unnodal}\; ({\it i.e.} without $(-n)$-curves, $n\ge 2$)\;
Del Pezzo surface with a connected real part $\R\Sig$, and let
$D\subset\Sig$ be a real ample divisor. Consider a generic
configuration~$\omega$ of $c_1(\Sig)\cdot D-1$ distinct real points
of $\Sig$. The set $R(D, \omega)$ of real ({\it i.e.} complex
conjugation invariant) rational curves $C\in|D|$ passing through the
points of~$\omega$ is finite, and all these curves are nodal and
irreducible. Put
$$W(\Sig,D,
\omega)=\sum_{C\in R(D, \omega)}(-1)^{s(C)}\ ,$$ where $s(C)$ is the
number of solitary nodes of $C$ ({\it i.e.} real points, where a
local equation of the curve can be written over $\R$ in the form
$x^2+y^2=0$). By Welschinger's theorem \cite{W1}, the number
$W(\Sig,D, \omega)$ does not depend on the choice of a generic
configuration $\omega$, and hence we simply write $W(\Sig,D)$,
omitting the configuration in the notation of this
Welschinger invariant. 

In what follows we shall use a generalized definition of the
Welschinger sign of a curve. Namely, let $C$ be a real algebraic
curve on a smooth real algebraic surface $\Sig$, and let
$\overline\bp\subset\Sig$ be a conjugation invariant finite subset.
Assume that $C$ has no singular local branches ({\it i.e.} is an
immersed curve). Then we define the Welschinger sign
\begin{equation}W_{\Sig,\overline\bp}(C)=(-1)^{s(C,\Sig,\overline\bp)},\quad\text{where}\quad
s(C,\Sig,\overline\bp)=\sum_{z\in\Sing(C')}s(C',z)\
,\label{e34}\end{equation} $C'$ being the strict transform of $C$
under the blow up of $\Sig$ at $\overline\bp$, and $s(C',z)$ is the
number of solitary nodes in a local $\del$-const deformation of the
singular point $z$ of $C'$ into $\del(C',z)$ nodes, where $\del$
denotes the $\del$-invariant of singularity (i.e., the maximal
possible number of nodes in its deformation). It is evident that
$s(C',z)$ is correctly defined modulo $2$, and hence
$W_{\Sig,\overline\bp}(C)$ is well-defined.

\medskip

{\bf Organization of the material.} In section \ref{ptc}, we set
forth the geometry of plane tropical curves adapted to our purposes,
completing with the definition of weights of tropical curves which,
in the complex case, designate the number of algebraic curves
associated with the given tropical curves in the further
patchworking theorem, and, in the real case, the contribution of the
real algebraic curves in the associated set to the Welschinger
number. In section \ref{sec3}, we provide two patchworking theorems,
the complex and the real one, in which we explicitly construct
algebraic curves associated to the tropical curves under
consideration.

\medskip

{\bf Acknowledgement}. The author was supported by the grant no.
465/04 from the Israel Science Foundation, a grant from the Higher
Council for Scientific Cooperation between France and Israel, and a
grant from the Tel Aviv University. This work was completed during
the author's stay at the Centre Interfacultaire Bernoulli, \'Ecole
Polytechnique F\'ed\'erale da Lausanne and at Laboratoire Emile
Picard, Universit\'e Paul Sabatier, Toulouse. The author thanks the
CIB-EPFL and UPS for the hospitality and excellent working
conditions. Special thank are due to I. Itenberg, who pointed out a
mistake in the preliminary version of Theorem \ref{t2}. Finally, I
express my gratitude to the unknown referee for numerous remarks,
corrections, and suggestions.

\section{Parameterized plane tropical curves}\label{ptc}

For the reader's convenience, we remind here some basic definitions
and facts about tropical curves, which we shall use in the sequel.
The details can be found in \cite{Mi1,Mi,RST}.

\subsection{Definition} An {\bf abstract tropical curve} is a
compact graph $\oGam$ without divalent vertices and isolated points
such that $\Gam=\oGam\backslash \Gam^0_\infty$, where
$\Gam^0_\infty$ is the set of univalent vertices, is a metric graph
whose closed edges are isometric to closed segments in $\R$, and
non-closed edges {\bf $\Gam$-ends} are isometric to rays in $\R$ or
to $\R$ itself. Denote by $\oGam^0$, resp. $\Gam^0$, the set of
vertices of $\oGam$, resp. $\Gam$, and split the set $\oGam^1$ of
edges of $\oGam$ into $\Gam^1_\infty$, the set of the $\Gam$-ends,
and $\Gam^1$, the set of closed (finite length) edges of $\Gam$. The
{\bf genus} of $\Gam$ is $g=b_1(\Gam)-b_0(\Gam)+1$.

A {\bf plane parameterized tropical curve} (shortly {\bf PPT-curve})
is a pair $(\oGam,h)$, where $\oGam$ is an abstract tropical curve
and $h:\Gam\to\R^2$ is a continuous map whose restriction to any
edge of $\Gam$ is a non-zero $\Z$-affine map and which satisfies the
following {\bf balancing} and {\bf nondegeneracy} conditions at any
vertex $v$ of $\Gam$: For each $v\in\Gam^0$,
\begin{equation}\sum_{v\in e,\ e\in\oGam^1}dh_v(\tau_v(e))=0\
,\label{e1}\end{equation}and $$\Span\{dh_v(\tau_v(e))\ ;\ v\in e,\
e\in\oGam^1\}=\R^2\ ,$$ where $\tau_v(e)$ is the unit tangent vector
to an edge $e$ at the vertex $v$. The {\bf degree} of a PPT-curve
$(\oGam,h)$ is the unordered multi-set of vectors $\{dh(\tau(e))\ :\
e\in\Gam^1_\infty\}$, where $\tau(e)$ denotes the unit tangent
vector of a $\Gam$-end $e$ pointing to the univalent vertex.

Observe that
\begin{equation}\sum_{e\in\Gam^1_\infty}dh(\tau(e))=0\ ,\label{e8}\end{equation} what
immediately follows from the balancing condition (\ref{e1}). We
shall also use another form of the $\Gam$-end-balancing condition.
For each $\Gam$-end $e$ pick any point \mbox{$\bx_e\in
h(e\backslash\Gam^0_\infty)$}. Then
\begin{equation}\sum_{e\in\Gam^1_\infty}\langle R_{\pi/2}(dh(\tau(e))),\bx_e\rangle=0\
,\label{e10}\end{equation} where $R_{\pi/2}$ is the (positive)
rotation by $\pi/2$. This is an elementary consequence of the stuff
discussed in the next section: one can lift a PPT-curve to a plane
algebraic curve over a non-Archimedean field, consider the defining
polynomial, and then use the fact that the product of the roots of
the (quasihomogeneous) truncations of this polynomial on the sides
of its Newton polygon equals $1$. We leave details to the reader.

Since $dh_v((\tau_v(e))\in\Z^2$, we have a well-defined positive
weight function \mbox{$w:\oGam^1\to\Z$} in the relation
$dh_v(\tau_v(e))=w(e)\bu_v(e)$ with $\bu_v(e)$ being the primitive
integral tangent vector to $h(e)$, emanating from $h(v)$. In the
sequel, when modifying tropical curves we speak of changes of edge
weights, which in terms of $h$ and $\oGam$ means that $h$ remains
unchanged whereas the metric on the chosen edges is multiplied by a
constant.

Observe that a connected component of $\oGam\backslash F$, where $F$
is finite, naturally induces a new PPT-curve (further on referred to
as {\bf induced}) when making the metric on the non-closed edges of
that component complete and respectively correcting the map $h$ on
these edges. These induced curves and the unions of few of them,
coming from the same $\oGam\backslash F$, are called PPT-curves {\bf
subordinate} to $(\oGam,h)$.

The {\bf deformation space} ${\cal M}(\oGam,h)$ of a PPT-curve
$(\oGam,h)$ is obtained by variation of the length of the finite
edges of $\Gam$ and combining $h$ with shifts. It can be identified
with an open rational convex polyhedron in Euclidean space, and its
closure $\overline{\cal M}(\oGam,h)$ can be obtained by adding the
boundary of that polyhedron which corresponds to PPT-curves with
some edges $e\in\Gam^1$ contracted into points.

Deformation equivalent PPT-curves are often called to be of the same
{\bf combinatorial type}. The degree and the genus are invariants of
the combinatorial type as well as the following characteristics. We
call a PPT-curve $(\oGam,h)$
\begin{itemize}\item {\bf irreducible}\; if $\Gam$ is connected, \item {\bf
simple}\; if $\Gam$ is trivalent, and \item {\bf pseudo-simple}\;
if, for any vertex $v\in\Gam^0$ incident to $m>3$ edges
$e_1,e_2,...,e_m$, one has $\bu_v(e_1)\ne \bu_v(e_j)$, $1<j\le m$,
and only two distinct vectors among
$\bu_v(e_2),...,\bu_v(e_m)$.\end{itemize} In the latter case, an
edge $e_i$ emanating from a vertex $v\in\Gam^0$ of valency $m>3$ is
called {\bf simple}, if $\bu_v(e_i)\ne\bu_v(e_j)$ for all $j\ne i$,
and is called {\bf multiple} otherwise.

\subsection{Newton polygon and its subdivision dual to a plane tropical curve}\label{sec4}
Given a PPT-curve $Q=(\Gam,h)$, the image $T=h(\Gam)\subset\R^2$ is
a finite planar graph, which supports an embedded plane tropical
curve (shortly {\bf EPT}-curve) $h_*Q:=(T,h_*w)$ with the (edge)
weight function
$$h_*w:T^1\to\Z,\quad h_*w(E)=\sum_{e\in\oGam^1,\ h(e)\supset E}
w(e)\ .$$ The respective balancing condition immediately follows
from (\ref{e1}). Furthermore, there exists a convex lattice
polygon $\Del\subset\R^2$ (different from a point) and a convex
piece-wise linear function
\begin{equation}f_T:\R^2\to\R,\quad
f(\bx)=\max_{\omega\in\Del\cap\Z^2}(\langle
\omega,\bx\rangle+c_\omega),\quad \bx\in\R^2\
,\label{e2}\end{equation} such that
\begin{itemize}\item $T$ is the corner locus of $f_T$, \item for
any two linearity domains $D_1,D_2$ of $f_T$, corresponding to
linear functions in formula (\ref{e2}) with gradients
$\omega_1,\omega_2$, respectively, and having a common edge
$E=D_1\cap D_2$ of $T$, it hods
$\omega_2-\omega_1=h_*w(E)\cdot\bu(E)$, where $\bu(E)$ is the
primitive integral vector orthogonal to $E$ and directed from $D_1$
to $D_2$.\end{itemize} Here the polygon $\Del$, called the {\bf
Newton polygon}\; of $Q$, is defined uniquely up to a shift in
$\R^2$, and $f_T$ is defined uniquely up to addition of a linear
affine function.

The Legendre dual to $f_T$ function $\nu_T:\Delta\to\R$ is convex
piece-wise linear, and its linearity domains define a subdivision
$S_T$ of $\Del$ into convex lattice subpolygons. This subdivision
$S_T$ is dual to the pair $(\R^2,T)$ in the following way: there is
1-to-1 correspondence between the faces of subdivision of $\R^2$
determined by $T$ and the faces of subdivision $S_T$ such that (i)
the sum of the dimensions of dual faces is $2$, (ii) the
correspondence inverts the incidence relation, (iii) the dual edges
of $T$ and of $S_T$ are orthogonal, and the weight of an edge of $T$
equals the lattice length of the dual edge of $S_T$. In particular,
If $V=(\alp,\bet)$ is a vertex of $T$, then
$\nabla\nu_T=(-\alp,-\bet)$ along the dual polygon $\Del_V$ of the
subdivision $S_T$.

Furthermore, we can obtain an extra information on the subdivision
$S_T$ out of the original PPT-curve $Q$. Namely,
\begin{itemize}\item with each edge $e\in\oGam^1$ we associate
a lattice segment $\sig_e$ which is orthogonal to $h(e)$ and
satisfies $|\sig_e|=w(e)$,
\item with each vertex $v\in\Gam^0$ we associate
a convex lattice polygon $\Del_v$, whose sides are suitable
translates of the segments $\sig_e$, $e\in\oGam^1$, $v\in e$.
Denote by $\sig_{v,e}$ the side of $\Del_v$, which is a translate
of $\sig_e$ and whose outward normal is $dh_v(\tau_v(e))$.
\end{itemize} Let a polygon $\Del_V$ of the subdivision $S_T$ be dual to a
vertex $V$ of $T$. Then (up to a shift)
\begin{equation}\Del_V=\sum_{\renewcommand{\arraystretch}{0.6}
\begin{array}{cc}
\scriptstyle{e\in\oGam^1}\\
\scriptstyle{\Int(e)\cap h^{-1}(V)\ne\emptyset}
\end{array}}\sig_e+\sum_{\renewcommand{\arraystretch}{0.6}
\begin{array}{cc}
\scriptstyle{v\in\Gam^0}\\
\scriptstyle{h(v)=V}
\end{array}}\Del_v\ .\label{e3}\end{equation} In this connection, we
can speak on $\nabla\nu_T$ along the polygons $\Del_v$ appearing in
(\ref{e3}).

A EPT curve $T$ is called {\bf nodal}, if the dual subdivision $S_T$
consists of triangles and parallelograms, i.e., when the
non-trivalent vertices of $T$ are locally intersections of two
straight lines. A nodal EPT curve canonically lifts into a simple
PPT curve when one resolves all nodes of the given curve.

\subsection{Compactified tropical curves}

For a given convex lattice polygon $\Del$, different from a point,
we define a compactification $\R^2_\Del$ of $\R^2$ in the following
way. If $\dim\Del=2$, we identify $\R^2$ with the positive orthant
$(\R_{>0})^2$ by the coordinate-wise exponentiation, then identify
$(\R_{>0})^2$ with the interior of
$\R^2_\Del:=\Tor_\R(\Del)_+\simeq\Del$, the non-negative part of the
real toric variety $\Tor_\R(\Del)$, via the moment map
$$\mu(x)=\frac{\sum_{\omega\in\Del\cap\Z^2}x^\omega\omega}{\sum_{\omega\in
\Del\cap\Z^2}x^\omega},\quad x\in(\R_{>0})^2\ .$$ If $\Del$ is a
segment, then we take $\Del'=\Del\times\sig$, $\sig$ being a
transverse lattice segment, and define $\R^2_\Del$ as the quotient
of $\R^2_{\Del'}$ by contracting the sides parallel to $\sig$. We
observe that the rays in $\R^2$, directed by an external normal
$\bu$ to a side $\sig$ of $\Del$ and emanating from distinct points
on a line, transverse to $\sig$, close up at distinct points on the
part of $\partial(\R^2_\Del)$, corresponding to the interior of
$\sig$ in the above construction\footnote{Clearly, the rays directed
by vectors distinct from any exterior normal to sides on $\Del$
close up at respective vertices of $\R^2_\Del$.}.

So, we can naturally compactify a PPT-curve $(\Gam,h)$ into
$(\oGam,\oh)$, by extending $h$ up to a map
$\oh:\oGam\to\R^2_\Del$.

\subsection{Marked tropical curves} An abstract
tropical curve with $n$ marked points is a pair $(\oGam,G)$, where
$\oGam$ is an abstract tropical curve and $G=(\gam_1,...,\gam_n)$ is
an ordered $n$-tuple of distinct points of $\oGam$. We say that a
marked tropical curve $(\oGam,G)$ is {\bf regular}\; if each
connected component of $\oGam\backslash G$ is a tree containing
precisely one vertex from $\Gam^0_\infty$. Furthermore, a marked
tropical curve $(\oGam,G)$ is called \begin{itemize}\item {\bf
end-marked}, if $G\cap\Gam^0=\emptyset$ and the points of $G$ lie on
the ends of $\oGam$, one on each end, \item {\bf regularly
end-marked}, if $G\cap\Gam^0=\emptyset$, the points of $G$ lie on
the ends of $\oGam$, and $(\oGam,G)$ is regular.\end{itemize}

A parametrization of a (compact) plane tropical curves with marked
points is a triple $(\oGam,G,\oh)$, where $(\oGam,G)$ is a marked
abstract tropical curve, and $(\oGam,\oh)$ is a PPT-curve. We define
the deformation space ${\cal M}(\oGam,G,\oh)\subset{\cal
M}(\oGam,\oh)$ by fixing the combinatorial type of the pair
$(\oGam,G)$ ($G$ being an ordered sequence). It can be identified
with a convex polyhedron in $\R^N$, where the coordinates designate
the two coordinates of the image $\oh(v)$ of a fixed vertex
$v\in\Gam^0$, the lengths of the edges $e\in\Gam^1$, and the
distances between the marked points lying inside edges of $\Gam$ to
some fixed points inside these edges (chosen one on each edge), cf.
\cite{GM1}. Further on the deformation type of a marked PPT-curve is
called a {\bf combinatorial} type.

\begin{lemma}\label{l1}
Let $\Del$ be a convex lattice polygon, $X=(\bx_1,...,\bx_n)$ a
sequence of points in $\R^2_\Del$ (not necessarily distinct). Then
there exists at most one $n$-marked regular PPT-curve
$(\oGam,G,\oh)$ with the Newton polygon $\Del$ and with a fixed
combinatorial type, such that $\oh(\gam_i)=\bx_i$, $\gam_i\in G$,
$i=1,...,n$.
\end{lemma}

{\bf Proof.} If such a marked PPT-curve exists, it is sufficient to
uniquely restore each connected component of $\oGam\backslash G$,
and hence, the general situation reduces to the case of an
irreducible rational PPT-curve (a subordinate curve defined by such
a connected component) with $|\Gam^0_\infty|-1=|\Gam^1_\infty|-1$
marked univalent vertices. We proceed by induction on $|\oGam^1|$.
The base of induction, {\it i.e.} the case $|\oGam^1|=1$ is evident.
Assume that $|\oGam^1|>1$.

If there are two $\Gam$-ends $e_1,e_2$ with marked points that
emanate from one vertex $v\in\Gam^0$ and are mapped into the same
straight line by $\oh$, then either $\oh(e_1)=\oh(e_2)$, in which
case we replace $e_1,e_2$ by one end of weight $w(e_1)+w(e_2)$ and
respectively replace two marked points by one, thus, reducing
$|\oGam^1|$ by $1$ and keeping the irreducibility and the
rationality of the tropical curve, or $\oh(e_1)$ and $\oh(e_2)$ are
the opposite rays emanating from $\oh(v)$, in which case we remove
the $\Gam$-end with lesser weight, leaving the other with weight
$|w(e_1)-w(e_2)|$, thus, reducing $|\oGam^1|$ by $1$ or $2$.

If there are no $\Gam$-ends as above, from
$$|\Gam^0|-|\Gam^1|=1\quad\text{and}\quad 3\cdot|\Gam^0|\le
2\cdot|\Gam^1|+|\Gam^1_\infty|$$ we deduce that
$|\Gam^0|\le|\Gam^1_\infty|-2$. Hence there are two non-parallel
$\Gam$-ends with marked points which merge to a common vertex
$v\in\Gam^0$, which thereby is determined uniquely. So, we remove
the above $\Gam$-ends and the vertex $v$ from $\oGam$, then extend
the other edges of $\oGam$ coming to $v$ up to new ends and mark on
them the points mapped to $h(v)$. Thus, the induction assumption
completes the proof. \proofend

\subsection{Tropically generic configurations of points} Let $\Del$ be a convex lattice polygon,
$\overline\bx=(\bx_1,...,\bx_k)$ a sequence of distinct points in
$\R^2_\Del$ such that $\bx_i\in\sig_i$, $1\le i\le r$, where
$\sig_1,...,\sig_r\subset\R^2_\Del$ correspond to certain sides of
$\Del$, and $\bx_i\in\R^2\subset\R^2_\Del$, $r<i\le k$. Let
$\overline m=(m_1,...,m_k)$ be a sequence of non-negative
integers, called weights of the points $\bx_1,...,\bx_k$,
respectively. A subconfiguration of $(\overline\bx,\overline m)$
is a configuration $(\overline\bx,\overline m')$ with $\overline
m'\le\overline m$ (component-wise).

Let ${\cal C}$ be a combinatorial type of an irreducible end-marked
PPT-curve with Newton polygon $\Del$, with $m=m_1+...+m_k$
$\Gam$-ends and marked points $\gam_1,...,\gam_m$. A weighted
configuration $(\overline\bx,\overline m)$ is called ${\cal C}$-{\bf
generic}, if there is no end-marked irreducible PPT-curve
$(\oGam,G,\oh)$ of type ${\cal C}$ such that
$\oh(G)=(\overline\bx,\overline m)$, {\it i.e.}
$$\oh(\gam_i)=\bx_i,\quad\sum_{j<i}m_j<i\le\sum_{j\le i}m_j,\quad
i=1,...,k\ .$$ A weighted configuration $(\overline\bx,\overline m)$
is called $\Del$-{\bf generic}, if it together with all its
subconfigurations is generic with respect to the combinatorial types
of end-marked irreducible PPT-curves which have
$m\le|\partial\Del\cap\Z^2|$ $\Gam$-ends and directing vectors of
all edges orthogonal to integral segments in $\Del$. A
(non-weighted) configuration $\overline\bx$ is called $\Del$-{\bf
generic}\; if all possible weighted configurations
$(\overline\bx,\overline m)$ are $\Del$-generic.

\begin{lemma}\label{l3}
The $\Del$-generic configurations with rational coordinates
\mbox{$\overline\bx=(\bx_1,...,\bx_k)\subset\R^2_\Del$} such that
$\bx_i\in\sig_i$, $1\le i\le r$, $\bx_i\in\R^2$, $r<i\le k$, form
a dense subset of $\sig_1\times...\times\sig_r\times(\R^2)^{k-r}$.
\end{lemma}

{\bf Proof}. Notice that there are only finitely many (up to the
choice of edge weights) combinatorial types of end-marked
irreducible PPT-curves under consideration and only finitely many
weight collections $\overline m$ to consider. We shall prove that,
for any such combinatorial type ${\cal C}$ of end-marked irreducible
PPT-curves, the image of the natural evaluation map $\Ev:{\cal
M}({\cal C})\to\sig_1\times...\times\sig_r\times(\R^2)^{k-r}$ is
nowhere dense, and hence is a finite polyhedral complex of a
positive codimension. This would suffice for the proof of Lemma due
to the aforementioned finiteness.

Thus, assuming that an end-marked irreducible curve
$(\overline\Gam,G,\overline h)$ of type ${\cal C}$ matches a
weighted rational configuration $(\overline\bx,\overline m)$, we
shall show that this imposes a nontrivial relation on the
coordinates of the points of $\overline\bx$, and hence complete the
proof. Clearly, any point $\bx_i\in\overline\bx$ lying on
$\partial\R^2_\Del$ ($\Del$ being the Newton polygon of
$(\overline\Gam,G,\overline h)$) is a univalent vertex for some ends
of $\overline\Gam$, whose $\overline h$-images lie on the same
straight line. So, pushing all the points of
$\overline\bx\cap\partial\R^2_\Del$ along the corresponding lines,
we can make $\overline\bx\subset\R^2$. Take an irrational vector
$a\in\R^2$ and pick the point $\bx_i\in\overline\bx$ with the
maximal value of the functional $\langle a,\bx\rangle$. Notice that
there are no vertex $v\in\Gam^0$ with $\langle
a,h(v)\rangle\ge\langle a,\bx_i\rangle$, since otherwise, due to the
balancing condition (\ref{e1}), one would find an end $e$ of $\Gam$
with $h(e)$ lying entirely in the half-plane $\langle
a,\bx\rangle>\langle a,\bx_i\rangle$ contrary the the assumptions
made. Hence, for each end $e\in\Gam^1_\infty$ with $h(e)$ passing
through $\bx_i$, we have $\langle a,\tau(e)\rangle>0$, which yields
that $$\sum_{\renewcommand{\arraystretch}{0.6}
\begin{array}{cc}
\scriptstyle{e\in\Gam^1_\infty}\\
\scriptstyle{\bx_i\in h(e)}
\end{array}}m_e\cdot dh(\tau(e))\ne 0$$ for any positive integers
$m_e$, and which finally implies that the coordinates of $\bx_i$
nontrivially enter relation (\ref{e10}). \proofend

\begin{lemma}\label{l2}
Let $\overline\bx$ be a $\Del$-generic configuration of points,
$Q=(\oGam,G,\oh)$ a marked regular PPT-curve with Newton polygon
$\Del$ which matches $\overline\bx$. Then
\begin{enumerate}\item[(i)] $(\oh)^{-1}(\overline\bx)=G$;\item[(ii)]
if $K$ is a connected component of $\oGam\backslash G$, then its
edges can be oriented so that

- the edges merging to marked points, emanate from these points,

- the unmarked $\Gam$-end is oriented towards its univalent
endpoint,

- from any vertex $v\in K^0$ emanates precisely one edge, and this
edge is simple.\end{enumerate}
\end{lemma}

\begin{remark}\label{r1} It follows from Lemma \ref{l2} that if an edge of $\Gam$ is
multiple for both of its endpoints and contains a marked point
inside which matches a point $\bx\in\overline\bx$, then all the
other edges joining the same vertices contain marked points matching
$\bx$.
\end{remark}

{\bf Proof of Lemma \ref{l2}}. (i) Assume that there is a point
$\gam\in(\oh)^{-1}(\overline\bx)\backslash G$. It belongs to a
component $K$ of $\oGam\backslash G$, which is a tree due to the
regularity of the considered marked tropical curve, and hence, is
cut by $\gam$ into two trees $K_1,K_2$, and only one of them, say
$K_1$ contains a $\Gam$-end free of marked points. Then, marking the
new point $\gam$, we obtain that the irreducible (rational)
PPT-curve induced by $K_2$ is end-marked and matches a
subconfiguration of $\overline\bx$ contrary to its
$\Del$-genericity.

(ii) Observe that the image of the unmarked ray does not coincide
with the image of any other edge of $K$, what immediately follows
from the statement (i).

Next we notice that, if $p$ is a vertex of $K$, $e$ a multiple edge
merging to $p$, then the connected component $K(e)$ of
$K\backslash\{p\}$, starting with the edge $e$ does not contain the
unmarked $K$-end. Indeed, otherwise, we consider another edge $e'$
of $K$ merging to $p$ so that $\bu_p(e')=\bu_p(e)$. Then we take the
graph $K\backslash K(e)$, which after a suitable modification of the
weights of the edges of $K(e')$ (the component of $K\backslash\{p\}$
starting with $e'$) induces an end-marked (rational) PPT-curve
matching the $\Del$-generic configuration $\overline\bx$, thus, a
contradiction.

It follows from the latter observation that $K$ has no edge being
multiple for both of its endpoints. Indeed, otherwise we would have
two vertices $v_1,v_2\in K^0$, joined by an edge $e\in K^1$,
multiple for both $v_1$ and $v_2$, and then would obtain that the
unmarked $K$-end is contained either in the component of
$K\backslash\{v_1\}$ starting with $e$, or in the component of
$K\backslash\{v_2\}$ starting with $e$, contrary to the above
conclusion.

Finally, we define an orientation of the edges of $K$, opposite to
the required one. Start with the unmarked $K$-end and orient it
towards its multivalent endpoint. In any other step, coming to a
vertex $v\in K^0$ along some edge, we orient all other edges merging
to $v$ outwards. Since $K$ is a tree, the orientation smoothly
extends to all of its edges. The preceding observations confirm that
any edge $e$ oriented in this manner towards a vertex $v\in E$ is
simple for $v$. \proofend

\subsection{Weights of marked pseudo-simple regular
PPT-curves}\label{sec1} In this section, $Q=(\oGam,G,\oh)$ is always
a regular marked pseudo-simple PPT-curve. Denote
$G_\infty=G\cap\Gam^0_\infty$ and $G_0=G\backslash G_\infty$, and
put $\overline\bx=\oh(G)$, $\overline\bx^\infty=\oh(G_\infty)$.
Throughout this section we assume that
{\it\begin{enumerate}\item[(T1)] no edge of $\Gam$ is multiple for
two vertices of $\Gam$,
\item[(T2)] $G_0$ does not contain vertices of valency $>3$,
\item[(T3)] $\overline\bx$ is $\Del$-generic.\end{enumerate}}

In particular, by Lemma \ref{l2}, we have that
$(\oh)^{-1}(\overline\bx)=G$.

{\bf Complex weights.} We define the {\bf complex weight} of a
PTT-curve $Q=(\oGam,G,\oh)$ as
\begin{equation}M(Q)=\prod_{v\in\Gam^0}M(Q,v)\cdot\prod_{e\in\Gam^1}M(Q,e)\cdot\prod_{\gam\in
G}M(Q,\gam)\ ,\label{e15}\end{equation} where the values
$M(Q,v),M(Q,e),M(Q,\gam)$ are computed along the following rules.

(M1) $M(Q,e)=w(e)$ for each edge $e\in\Gam^1$.

(M2) $M(Q,\gam)=1$ for each $\gam\in
G\cap(\Gam^0\cup\Gam^0_\infty)$, and $M(Q,\gam)=w(e)$ for each
$\gam\in G\backslash(\Gam^0\cup\Gam^0_\infty)$, $\gam\in
e\in\oGam^1$.

To define $M(Q,v)$, $v\in\Gam^0$, introduce some notation: denote
by $\Del_v$ the lattice triangle whose boundary is combined of the
vectors $dh(\tau_v(e))$, rotated clockwise by $\pi/2$, where $e$
runs over all the edges of $\oGam$ emanating from $v$. Next, we
put:

(M3) If $v\in\Gam^0\cap G$, then $M(Q,v)=|\Del_v|$.

(M4) If $v\in\Gam^0\backslash G$ is trivalent, then it belongs to a
connected component $K$ of $\oGam\backslash G$ which we orient as in
Lemma \ref{l2}(ii) and thus define two edges $e_1,e_2\in\oGam^1$
merging to $v$. In this case we put
$M(Q,v)=|\Del_v|(w(e_1)w(e_2))^{-1}$.

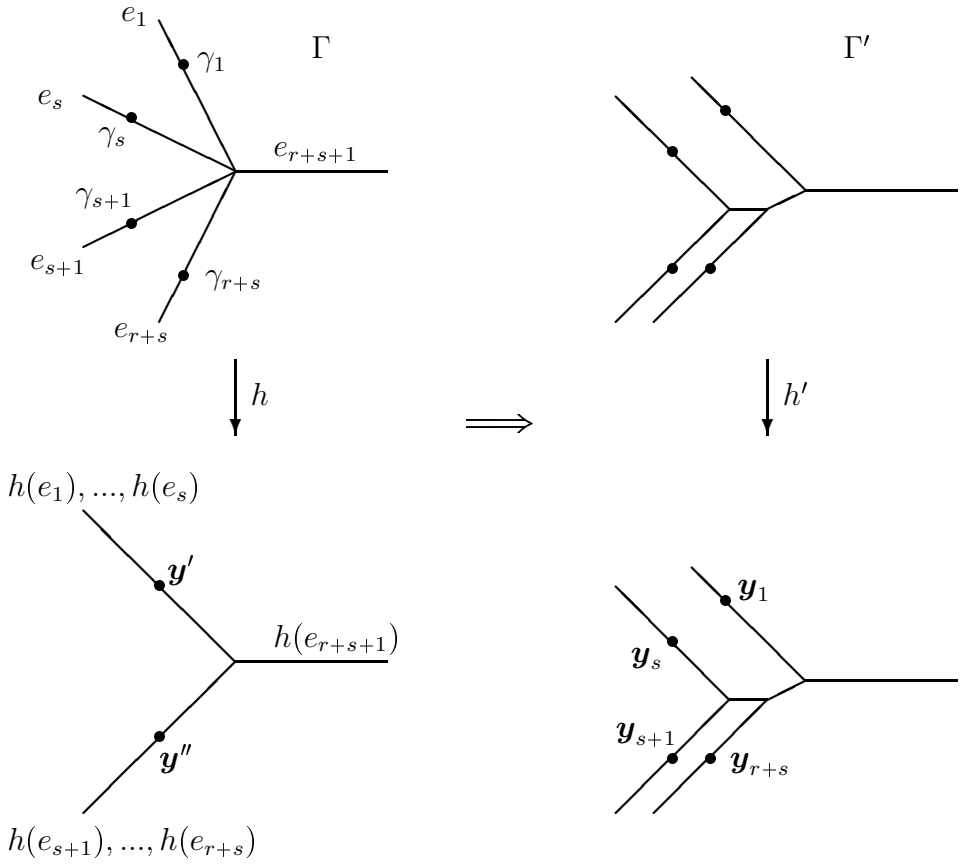
\begin{figure}
\setlength{\unitlength}{1cm}
\begin{picture}(13,12)(0,0)
\thicklines\put(1,1){\line(1,1){2}}\put(3,3){\line(-1,1){2}}
\put(3,3){\line(1,0){2}}\put(2,7.5){\line(1,2){1}}
\put(1,8.5){\line(2,1){2}}\put(3,9.5){\line(1,0){2}}\put(3,9.5){\line(-1,2){1}}\put(3,9.5){\line(-2,1){2}}
\put(8.5,1){\line(1,1){1.5}}\put(8,1){\line(1,1){1.5}}\put(9.5,2.5){\line(1,0){0.5}}\put(9.5,2.5){\line(-1,1){1.5}}
\put(10,2.5){\line(2,1){0.5}}\put(10.5,2.75){\line(-1,1){1.5}}\put(10.5,2.75){\line(1,0){2}}
\put(8.5,7.5){\line(1,1){1.5}}\put(8,7.5){\line(1,1){1.5}}\put(9.5,9){\line(1,0){0.5}}\put(9.5,9){\line(-1,1){1.5}}
\put(10,9){\line(2,1){0.5}}\put(10.5,9.25){\line(-1,1){1.5}}\put(10.5,9.25){\line(1,0){2}}
\put(3,7){\vector(0,-1){1}}\put(3.2,6.4){$h$}\put(10,7){\vector(0,-1){1}}\put(10.2,6.4){$h'$}
\put(6,6){\Large{$\Longrightarrow$}}\put(4,11){$\Gam$}\put(11,11){$\Gam'$}
\put(1.9,1.9){$\bullet$}\put(1.9,3.9){$\bullet$}\put(2.22,8){$\bullet$}\put(1.53,8.7){$\bullet$}
\put(2.22,10.8){$\bullet$}\put(1.53,10.1){$\bullet$}\put(9.15,1.6){$\bullet$}\put(8.65,1.6){$\bullet$}
\put(8.65,3.15){$\bullet$}\put(9.35,3.7){$\bullet$}\put(9.15,8.1){$\bullet$}\put(8.65,8.1){$\bullet$}
\put(8.65,9.65){$\bullet$}\put(9.35,10.2){$\bullet$}\put(2,1.6){$\by''$}\put(2.1,4.1){$\by'$}
\put(3.5,3.2){$h(e_{r+s+1})$}\put(0,0.5){$h(e_{s+1}),...,h(e_{r+s})$}
\put(0,5.2){$h(e_1),...,h(e_s)$}\put(1.4,7.3){$e_{r+s}$}\put(0.3,8.2){$e_{s+1}$}
\put(0.4,10.4){$e_s$}\put(1.5,11.5){$e_1$}\put(3.5,9.7){$e_{r+s+1}$}
\put(2.6,8){$\gam_{r+s}$}\put(0.9,9.1){$\gam_{s+1}$}\put(1.2,9.9){$\gam_s$}\put(2.5,10.9){$\gam_1$}
\put(8,2){$\by_{s+1}$}\put(9.5,1.6){$\by_{r+s}$}\put(8.2,3){$\by_s$}\put(9.6,3.9){$\by_1$}
\end{picture}
\caption{Local deformation of a tropical curve in (M5)}\label{f1}
\end{figure}

(M5) Let $v\in\Gam^0$ be of valency $s+r+1>3$, where $1\le r\le s$,
$2\le s$, and let $e_i$, $i=1,...,s+r+1$, be all the edges with
endpoint $v$ so that the edges $e_i$, $1\le i\le s$, have a common
directing vector $\bu_v(e_1)$, the edges $e_i$, $s<i\le s+r$, have a
common directing vector $\bu_v(e_{s+1})$, and $e_{s+r+1}$ is a
simple edge emanating from $v$ along the orientation of Lemma
\ref{l2}(ii). Consider a rational PPT-curve $Q_v$ induced by the
graph $\oGam_v=\{v\}\cup\bigcup_{i=1}^{s++r+1}e_i\subset\oGam$, pick
auxiliary marked points $\gam_i\in e_i\backslash\{v\}$,
$i=1,...,s+r$, in such a way that $\oh(\gam_i)=\by'\in\R^2$ as $1\le
i\le s$, and $\oh(\gam_i)=\by''\in\R^2$ as $s<i\le s+r$. Then we
replace $\by'$ (resp. $\by''$) by a generic set of distinct points
$\by_1,...,\by_s$ close to $\by'$ (resp. distinct points
$\by_{s+1},...,\by_{s+r}$ close to $\by''$), and take rational
regularly end-marked PPT-curves of degree
$\{dh(\tau_v(e_i))\}_{i=1,...,s+r+1}$ matching the configuration
$\by_1,...,\by_{s+r}$ so that the $\oh$-image of the $\Gam$-end of
weight $w(e_i)$ with the directing vector $\bu_v(e_i)$ passes
through the point $\by_i$, $i=1,...,s+r$ (see Figure
\ref{f1}).\footnote{Notice that by construction there is a canonical
1-to-1 correspondence between the ends of $Q_v$ and the ends of any
of the curves obtained in the deformation.} By \cite[Corollaries
2.24 and 4.12]{Mi}, the set ${\cal T}$ of these PPT-curves is
finite, and they all are simple. Then put
\begin{equation}M(Q,v)=\sum_{Q'\in{\cal
T}}M(Q')\ ,\label{e9}\end{equation} where all terms $M(Q')$ are
computed by formula (\ref{e15}) and the rules (M1)\;-(M4).

\begin{remark}\label{r2}
(1) We point out that the right-hand side of (\ref{e9}) does not
depend on the choice of the configuration $(\by_i)_{i=1,...,s+r+1}$,
what follows from \cite[Theorem 4.8]{GM1} (observe that the degree
of the evaluation map as in \cite[Definition 4.6]{GM1} coincides
with the right-hand side of (\ref{e9}) in our situation). Slightly
modifying the Mikhalkin correspondence theorem (\cite[Theorem
1]{Mi}), one can deduce that $M(Q,v)$ as defined in (\ref{e9})
equals the number of complex rational curves $C$ on the toric
surface $\Tor(\Del_v)$ such that \begin{itemize}\item $C$ belongs to
the tautological linear system $|{\cal L}_{\Del_v}|$,
\item for each side $\sig$ of $\Del_v$, the intersection points of
$C$ with the toric divisor $\Tor(\sig)\subset\Tor(\Del_v)$ are in
1-to-1 correspondence with the $\Gam$-ends of the tropical curves
from ${\cal T}$, orthogonal to $\sig$, $C$ is nonsingular along
$\Tor(\sig)$, and the intersection multiplicities are respectively
equal to the weights of the above $\Gam$-ends, \item $C$ passes
through a generic configuration of $s+r+1$ points in $\Tor(\Del_v)$.
\end{itemize}
Furthermore, ${\cal T}$ consists of just one curve as $r=1$. Indeed,
its dual subdivision of the Newton triangle $\Del_v$ must be as
described above with the order of segments dual to the parallel
$\Gam$-ends, which is determined uniquely by the disposition of the
points $\by_1,...,\by_s$.

(2) If $Q$ is simple ({\it i.e.} all the vertices of $\Gam$ are
trivalent, then (\ref{e15}) gives
\begin{equation}M(Q)=\frac{\prod_{v\in\Gam^o}|\Del_v|}{\prod_{\gam\in
G_\infty,\ \gam\in e\in\Gam^1_\infty}w(e)}\ ,
\label{e17}\end{equation} which generalizes Mikhalkin's weight
introduced in \cite[Definitions 2.16 and 4.15]{Mi}, and coincides
with the multiplicity of a tropical curve from \cite{GM1}.
\end{remark}

{\bf Real weights.} A PPT-curves $Q=(\oGam,G,\oh)$ equipped with an
additional structure, a continuous involution
$\conj:(\oGam,G,\oh)\righttoleftarrow$\; and a subdivision $G=\Re
G\cup\Im G$ invariant with respect to $\conj$, is called {\bf real}.

Clearly, $\Im\oGam:=\oGam\backslash\Re\oGam$, where
$\Re\oGam=\Fix\;\conj\big|_{\oGam}$ consists of two disjoint subsets
$\Im\oGam'$, $\Im\oGam''$ interchanged by $\conj$.

Given a real PPT-curve $Q$, we can construct a (usual) PPT-curve
$Q/\conj=(\oGam/\conj,G/\conj,\oh/\conj)$. Notice that the weights
of the edges obtained here by identifying $\Im\oGam'$ and
$\Im\oGam''$ are even. Conversely, given a (usual) PPT-curve
$Q=(\oGam,G,\oh)$ and a set $I(\oGam^1)\subset\oGam^1$, which
includes only edges of even weight, we construct a real PPT-curve
$Q'=(\oGam',G',\oh')$ as follows: (i) put $K=\bigcup_{e\in
I(\oGam^1)}e$ and obtain the graph $\oGam'$ by gluing up $\oGam$
with another copy $K'$ of $K$ at the vertices of $\Gam$, common for
$K$ and the closure of $\oGam\backslash K$, (ii) the map $\oh$
coincides on $K$ and $K'$, whereas the weights of the doubled edges
are divided by $2$ in order to keep the balancing condition, (iii)
the points of $G\cap K$ are respectively doubled to $K'$. Finally,
define an involution $\conj$ on $Q'$ interchanging $K$ and $K'$, and
define a subdivision $G'=\Re G'\cup\Im G'$.

We shall consider only real PPT-curves with the following
properties: {\it \begin{enumerate}\item[(R1)] $\Re\oGam$ is nonempty
and has no one-point connected component,
\item[(R2)] $\Im\oGam$
has only uni- and trivalent vertices (if nonempty),
\item[(R3)] the marked points $G_0\cap\Im\oGam$ are not vertices
of $\Im\oGam$, \item[(R4)] $\Im G\backslash\Im\oGam$ is empty or
consists of some trivalent vertices of $\Gam$, \item[(R5)] the
closure of any component of $\Im\oGam\backslash G$ contains a point
from $\Im G$.
\end{enumerate}}

Observe that the closure of $\Im\oGam$ joins $\Re\oGam$ at vertices
of valency $>3$ (which are not in $G$ by condition (T2)).

The real weight of a real PPT-curve $Q$ is defined as
\begin{equation}W(Q)=(-1)^{\ell_1}2^{\ell_2}\cdot
\prod_{v\in\Gam^0}W(Q,v)\cdot\prod_{e\in\oGam^1}W(Q,e)\cdot\prod_{\gam\in
G_0}W(Q,\gam)\ ,\label{e16}\end{equation} $$\ell_1=\frac{|\Re
G\cap\Im\oGam|}{2},\quad \ell_2=\frac{|\Im
G\cap\Im\oGam|-b_0(\Im\oGam)}{2}\ ,$$ with $W(Q,v),W(Q,e),W(Q,\gam)$
computed along the following rules:

(W1) For an edge $e\subset\Re \oGam$, put $W(Q,e)=0$ or $1$
according as $w(e)$ is even or odd. For an edge $e\in\Gam^1$,
$e\subset\overline{\Im\Gam}$, put $W(Q,e)W(Q,c(e))=w(e)$. For an
edge $e\in\Gam^1_\infty$, $e\subset\Im\oGam$, put $W(Q,e)=1$.

(W2) For $\gam\in G_0\cap\Re\oGam\backslash\Gam^0$, put
$W(Q,\gam)=1$. For $\gam\in \Re G\cap\Gam^0$, put $W(Q,\gam)=1$. For
$\gam\in\Im G$, $\gam=v\in\Gam^0$, put $W(Q,\gam)=|\Del_v|$. For
$\gam\in G_0$, $\gam\in e\subset\Im\oGam$, put
$W(Q,\gam)M(Q,c(\gam))=w(e)$.

(W3) For a vertex $v\in\Gam^0\cap\Im\oGam$, put
$W(Q,v)W(Q,\conj(v))=(-1)^{|\partial\Del_v\cap\Z^2|}M(Q,v)$ (see
condition (M4) for the definition of $M(Q,v)$). For a trivalent
vertex $v\in\Gam^0\cap\Re\oGam$, put
\mbox{$W(Q,v)=(-1)^{|\Int(\Del_v)\cap\Z^2|}$}.

(W4) For a four-valent vertex $v\in\Gam^0$ incident to two simple
edges from $\Re\oGam$ and two multiple edges $e',e''$ from
$\Im\oGam$, put
$W(Q,v)=(-1)^{|\Int(\Del_v)\cap\Z^2|}|\Del_v|/(2w(e'))$.

(W5) Let $v\in\Gam^0$ be of valency $>3$ incident to
\begin{itemize}\item a simple edge $e_1\subset\Re\oGam$,\item edges
$e_i\subset\Re\oGam$, $1<i\le r_1+1$, and $e'_i\subset\Im\oGam'$,
$e''_i\subset\Im\oGam''$, $1\le i\le s_1$, for some nonnegative
$r_1,s_1$, all with the same directing vector $\bu'\ne\bu_v(e_1)$,
and \item edges $e_i\subset\Re\oGam$, $r_1+1<i\le r_1+r_2+1$, and
$e'_i\subset\Im\oGam'$, $e''_i\subset\Im\oGam''$, $s_1< i\le
s_1+s_2$, for some nonnegative $r_2,s_2$ such that $r_2+2s_2\ge 2$,
all with the same directing vector $\bu''\ne\bu_v(e_1),\bu'$.
\end{itemize} Take the real PPT-curve $Q_v$ induced by $v$ and the
edges emanating from $v$, correspondingly restrict on $Q_v$ the
involution $\conj$, and introduce a finite $\conj$-invariant set of
marked points $G_v$ picking up one point on each edge emanating from
$v$ but $e_1$. Consider the PPT-curve $Q_v/\conj$ and perform with
it the deformation procedure described in (M5) (cf. Figure \ref{f1})
getting a finite set of simple rational regularly end-marked
PPT-curves. Any curve $\widetilde Q=(\widetilde\oGam,\widetilde
G,\widetilde\oh)$ from this set, we turn into a real PPT-curve.
Namely, first, we include into the set $I(\widetilde\oGam^1)$ all
the $\Gam$-ends which correspond to the $\Gam/\conj$-ends of $Q_v$
from $\Im\oGam_v/\conj$. Then we maximally extend the set
$I(\widetilde\oGam^1)$ in the following inductive procedure: if two
edges $f_1,f_2\in I(\widetilde\oGam^1)$ merge to a vertex
$p\in\widetilde\Gam^0$, then the third edge $f_3$, emanating from
$p$ should be added to $I(\widetilde\oGam^1)$. Clearly, by
construction, the weights of the edges $e\in I(\widetilde\oGam^1)$
are even; hence we can make a real PPT-curve $Q'=(\oGam',G',\oh')$,
letting $\Re G'=G'\cap\Re\oGam'$, $\Im G'=G'\cap\Im\oGam'$. Denoting
the final set of real PPT-curves by ${\cal T}$ and observing that
their real weight $W(Q')$ can be computed along the above rules
(W1)\;-(W4), we define
$$W(Q,v)=\sum_{Q'\in{\cal T}}W(Q')\ .$$ The fact that the latter expression does not depend on
the choice of the perturbation of the points $\by',\by''$ (cf.
construction in (M5) and Figure \ref{f1}) follows from a more
general statement proven in \cite{Sh08}.

\begin{remark}
(1) If $\conj=\Id$, $\Im G=\emptyset$, and $Q$ is simple, we obtain
the well-known formula: $W(Q)=0$ when $\oGam$ contains an even
weight edge, and $W(Q)=(-1)^a$,
$a=\sum_{v\in\Gam^0}|\Int(\Del_v)\cap\Z^2|$, when all the edge
weights of $\oGam$ are odd (cf. \cite[Definition 7.19]{Mi} or
\cite[Proposition 6.1]{Sh0}, where, in addition, $\deg Q$ consists
of only primitive integral vectors).

(2) If $Q$ is rational, $Q/\conj$ is simple and $G_\infty=\Re
G\cap\Gam^0=\Re G\cap\Im\oGam=\emptyset$, we obtain a generalization
of \cite[Formula (2.12)]{Sh1} (in version at arXiv:math/0406099).
Indeed, if $\Re\oGam$ contains an edge of even weight, we obtain
$W(Q)=0$ in (\ref{e16}) due to (W1), and accordingly we obtain
$w(Q/\conj)=0$ in \cite[Section 2.5]{Sh1} (in the notations
therein). If $\Re\oGam$ contains only edges of odd weight, then
\cite[Formula (2.12)]{Sh1} reads
\begin{equation}w(Q/\conj)=(-1)^{a+b}\prod_{v\in\Gam^0\cap\Im G}
|\Del_v|\cdot\prod_{v\in(\Gam/\conj)^0\cap\overline{(\Im\Gam/\conj)}}
\frac{|\Del_v|}{2}\label{e16n}\end{equation}
with $a=\sum_{v\in(\Gam/\conj)^0}|\Int(\Del_v)\cap\Z^2|$,
$b=|(\Gam/\conj)^0\cap(\Im\Gam/\conj)|$, whereas in (\ref{e16}) we
obtain $\ell_1=0$ by the assumption $\Re G\cap\Im\oGam=\emptyset$,
$\ell_2=|(\Gam/\conj)^0\cap(\Im\Gam/\conj)|$ due to the rationality
of $Q$ and simplicity of $Q/c$, and, furthermore, taking into
account that $w(e)=2w(e')=2w(\conj(e'))$ for
$e=(e'\cup\conj(e'))/\conj\in(\Gam\conj)^1$, $e'\in\Gam^1$,
$e'\subset\Im\Gam$, we compute the other factors in (\ref{e16}):
$$\prod_{v\in\Gam^0}W(Q,v)=\prod_{v\in\Gam^0\backslash\overline{\Im\Gam}}(-1)^{|\Int(\Del_v)\cap\Z^2|}
\cdot\prod_{\{v,\conj(v)\}\in(\Gam/\conj)^0\cap(\Im\Gam/\conj)}(-1)^{|\partial\Del_v\cap\Z^2|}M(Q,v)$$
$$\times\prod_{\renewcommand{\arraystretch}{0.6}
\begin{array}{cc}
\scriptstyle{v\in\Re\Gam\cap\overline{\Im\Gam/\conj}}\\
\scriptstyle{v\in e\subset\overline{\Im\Gam/c},\ e\in(\Gam/\conj)^1}
\end{array}}(-1)^{|\Int(\Del_v)\cap\Z^2|}\;\frac{|\Del_v|}{w(e)}$$
$$=\prod_{v\in(\Gam/\conj)^0}(-1)^{|\Int(\Del_v)\cap\Z^2|}\cdot(-4)^{|(\Gam/\conj)^0\cap
\Im\Gam/\conj|}\cdot
\prod_{v\in(\Gam/\conj)^0\cap\overline{\Im\Gam/\conj}}|\Del_v|$$
$$\times 2^{-|\Re\Gam\cap\overline{\Im\Gam}|}
\prod_{e\in(\Gam/\conj)^1,\
e\subset\overline{\Im\Gam/\conj}}\;\frac{2}{w(e)}\prod_{\renewcommand{\arraystretch}{0.6}
\begin{array}{cc}
\scriptstyle{e\in(\Gam/\conj)^1,\ e\subset\overline{\Im\Gam/\conj}}\\
\scriptstyle{e\cap G/\conj\ne\emptyset}
\end{array}}\frac{2}{w(e)}\ ,$$
$$\prod_{e\in\oGam^1}W(Q,e)=\prod_{e\in(\Gam/\conj)^1,\
e\subset\overline{\Im\Gam/\conj}}\;\frac{w(e)}{2}\prod_{\renewcommand{\arraystretch}{0.6}
\begin{array}{cc}
\scriptstyle{e\in(\Gam/\conj)^1,\ e\subset\overline{\Im\Gam/\conj}}\\
\scriptstyle{e\cap G/\conj\ne\emptyset}
\end{array}}\frac{w(e)}{2}\ ,$$ $$\prod_{\gam\in
G_0}W(Q,\gam)=\prod_{v\in\Im
G\cap\Re\Gam}|\Del_v|\prod_{\renewcommand{\arraystretch}{0.6}
\begin{array}{cc}
\scriptstyle{e\in(\Gam/\conj)^1,\ e\subset\overline{\Im\Gam/\conj}}\\
\scriptstyle{e\cap G/\conj\ne\emptyset}
\end{array}}\frac{w(e)}{2}\ ,$$ which altogether gives (with $a,b$ from (\ref{e16n})) $$W(Q)=(-1)^{a+b}
\prod_{v\in\Gam^0\cap\Im
G}|\Del_v|\cdot\prod_{v\in(\Gam/\conj)^0\cap\overline{(\Im\Gam/\conj)}}
\frac{|\Del_v|}{2}=w(Q/c)\ .$$
\end{remark}

\section{Patchworking theorem}\label{sec3}

\subsection{Patchworking data}\label{sec2}
{\bf Combinatorial-geometric part.} In the notation of section
\ref{sec1}, let $Q=(\oGam,G,\oh)$ be a pseudo-simple irreducible
regular marked PPT-curve of genus $g$ which has a nondegenerate
Newton polygon $\Del$ and which satisfies condition (T1)-(T3) of
section \ref{sec1}.

Let $G_0$ 
split into disjoint subsets $G_0=G^{(\mr)}_0\cup G^{(\dmp)}_0$ such
that $G^{(\mr)}_0\cap\Gam^0=\emptyset$ and $h(G_0^{(\mr)})\cap
h(G_0^{(\dmp)})=\emptyset$.
We equip
the points of $G_0$ with the following multiplicities:
\begin{itemize}\item if $\gam\in G^{(\mr)}$, put $\mt(\gam)=1$, \item
if $\gam\in G^{(\dmp)}$ is a (trivalent) vertex of $\Gam$, put
$\mt(\gam)=(1,1)$, \item if $\gam\in G^{(\dmp)}$ is not a vertex of
$\Gam$, put $\mt(\gam)=(1,0)$ or $(0,1)$.
\end{itemize}

In the sequel, by $\hat Q$ we denote the PPT-curve $Q$ equipped with
the subdivision $G_0=G^{(\mr)}_0\cup G^{(\dmp)}_0$ and the
multiplicity function $\mt(\gam)$, $\gam\in G_0$ as above.

\begin{definition}\label{newd}
A pair $\gam,\gam'$ of distinct points in $G_0$ is called {\bf
special} if $h(\gam)=h(\gam')$ and $\mt(\gam)=\mt(\gam')$. A pair of
parallel multiple edges $e,e'\in\oGam^1$ emanating from a vertex
$v\in\Gam^0$ of valency $>3$ is called {\bf special} if there are
disjoint open connected subsets $K,K'$ of $\Gam\backslash\{v\}$ such
that

- $K$ contains the germ of $e$ at $v$ and the point $\gam$, $K'$
contains the germ of $e'$ at $v$ and the point $\gam'$,

- there is a homeomorphism $\varphi:K\to K'$ satisfying
$h\big|_K=h\big|_{K'}\circ\varphi$.

A vertex $v\in\Gam^0$ incident to a special pair of edges is called
{\bf special}.
\end{definition}

\begin{figure}
\setlength{\unitlength}{1cm}
\begin{picture}(12,10)(0,0)
\thicklines\put(0.5,1.5){\line(1,1){1}}\put(1.5,2.5){\line(-1,1){1}}
\put(1.5,2.5){\line(1,0){3}}\put(3.5,2.5){\line(1,1){1}}
\put(3.5,2.5){\line(0,-1){1}}\put(0.5,5.5){\line(1,1){1}}\put(1.5,6.5){\line(-1,1){1}}\put(1.5,6.5){\line(1,0){3}}
\put(1.5,6.5){\line(2,1){2}}\put(3.5,7.5){\line(1,0){1}}\put(3.5,7.5){\line(1,1){0.5}}\put(7,1.5){\line(1,1){1}}
\put(8,2.5){\line(-1,1){1}}\put(8,2.5){\line(1,0){2}}\put(10,2.5){\line(1,1){1}}
\put(10,2.5){\line(1,-1){1}}\put(7,5.5){\line(1,1){1}}\put(8,6.5){\line(-1,1){1}}\put(8,6.5){\line(1,0){2}}
\put(10,6.5){\line(1,1){1.5}}\put(10,6.5){\line(1,-1){1}}\put(8,6.5){\line(2,1){2}}
\put(10,7.5){\line(1,1){1}}\put(10.05,7.45){\line(1,-1){0.4}}\put(10.55,6.95){\line(1,-1){1}}
\put(2.5,5){\vector(0,-1){1}}\put(2.7,4.4){$h$}\put(9,5){\vector(0,-1){1}}\put(9.2,4.4){$h$}
\put(2.4,2.38){$\bullet$}\put(2.4,6.38){$\bullet$}\put(2.4,6.88){$\bullet$}\put(10.4,1.9){$\bullet$}
\put(10.4,2.9){$\bullet$}\put(10.4,5.9){$\bullet$}\put(10.9,6.4){$\bullet$}
\put(10.4,7.9){$\bullet$}\put(10.9,7.38){$\bullet$}\put(1.5,6.1){$v$}\put(8,6.1){$v$}
\put(2.4,6.1){$\gam_1$}\put(2.3,7.2){$\gam_2$}\put(3.4,6.1){$e_1$}\put(3.1,7){$e_2$}
\put(9,6.1){$e_1$}\put(8.8,7.2){$e_2$}\put(10,5.7){$\gam_1$}\put(11.1,6.6){$\gam_2$}
\put(11.1,7.3){$\gam_3$}\put(10,8.1){$\gam_4$}\put(1.5,9){$m(\gam_1)=m(\gam_2)=1$}
\put(6.5,9){$m(\gam_1)=m(\gam_2)=m(\gam_3)=(1,0)$}\put(6.5,8.5){$m(\gam_4)=(0,1)$}
\put(2.5,1){\rm (a)}\put(9,1){\rm (b)}\put(0,0){\rm special pair of
marked points $(\gam_1,\gam_2)$, special pair of edges $(e_1,e_2)$,
special vertex $v$}
\end{picture}
\caption{Illustration to Definition \ref{newd}}\label{f2}
\end{figure}

Then we assume the following: {\it\begin{enumerate}
\item[(T4)] The edges in special pairs have weight $1$, and at
least one of the simple edges emanating from a special vertex has
weight $1$.
\item[(T5)] Let $e$, $e'$ be a special pair of edges emanating from
a vertex $v\in\Gam^0$, and let $K$, $K'$ be disjoint connected
subsets of $\Gam\backslash\{v\}$ as in Definition \ref{newd}; then
$K\cup K'$ contains at most one special pair of points of $G_0$.
\item[(T6)] A special pair of edges cannot be a pair of $\Gam$-ends 
and cannot be a pair of finite length edges
which end up at a special pair $\gam,\gam'\in\Gam^0$ such that
$h(\gam)=h(\gam')$ and $\mt(\gam)=\mt(\gam')=(1,1)$.
\item[(T7)] Let a vertex $v\in\Gam^0$ be a special vertex, and let
$\{e_1,...,e_s\}$ be a maximal (with respect to inclusion) set of
edges of $\Gam$ incident to $v$ and such that

- each edge $e_i$ contains a point $\gam_i\in G_0$, $1\le i\le s$,

- $h(\gam_1)=...=h(\gam_s)$ and $\mt(\gam_1)=...=\mt(\gam_s)$.

Suppose that $\dist(v,v_i)\le\dist(v,v_{i+1})$, $1\le i<s$,
$v_i\in\oGam^0\backslash\{v\}$ being the second vertex of $e_i$.
Then we require
\begin{equation}\dist(v,\gam_1)>\sum_{1\le i<s-1}\dist(\gam_i,v_i)+
2\cdot\dist(\gam_{s-1},v_{s-1})\ .\label{newe}\end{equation}
\end{enumerate}}

Notice that, in condition (T7), at most one edge $e_i$ is a
$\Gam$-end (cf. (T6)) and it must be $e_s$.

\medskip

Introduce also the semigroup
$$\Z_{\ge0}^\infty=\left\{\alp=(\alp_1,\alp_2,...)\ :\ \alp_i\in\Z,\ \alp_i\ge 0,\
i=1,2,...\ ,\ |\{i\ :\ \alp_i>0\}|<\infty\right\}\ ,$$ equipped with
two norms
$$\|\alp\|_0=\sum_{i=1}^\infty\alp_i,\quad\|\alp\|_1=\sum_{i=1}^\infty i\alp_i\ ,$$
and the partial order
$$\alp\ge\bet\quad\Leftrightarrow\quad\alp-\bet\in\Z_{\ge0}^\infty\ .$$ For each side $\sig$ of
$\Del$ introduce the vectors $\bet^\sig\in\Z_{\ge0}^\infty$ such
that the coordinate $\bet^\sig_i$ of $\bet^\sig$ equals the number
of the univalent vertices $v\in\Gam^0_\infty$ such that
$\oh(v)\in\sig\subset\R^2_\Del$ and $w(e)=i$ for the $\Gam$-end $e$
merging to $v$, for all $i=1,2,...$

\medskip
\noindent {\bf Algebraic part.} Let $\Sig=\Tor_\K(\Del)$. The
coordinate-wise valuation map $\val:(\K^*)^2\to\R^2$ naturally
extends up to $\val:\Sig\to\R^2_\Del$. Let
$\overline\bp\subset\Sig:=\Tor_\K(\Del)$ be finite and satisfy
$\val(\overline\bp)=\overline\bx=\oh(G)$.\footnote{This means, in
particular, that the points of $\overline\bx$ have rational
coordinates.} Suppose that {\it\begin{enumerate}\item[(A1)] each
point $\bx\in h(G^{(\mr)})\subset\overline\bx$ has a unique preimage
in $\overline\bp$, \item[(A2)] the preimage of each point $\bx\in
h(G^{(\dmp)})\subset\overline\bx$ consists of an ordered pair of
points $\bp_{1,\bx},\bp_{2,\bx}\in\overline\bp$, \item[(A3)] there
is a bijection $\psi:\overline\bp^\infty\to G_\infty$, where
$\overline\bp^\infty:=\val^{-1}(\overline\bx^\infty)$,
$\overline\bx^\infty=\overline h(G_\infty)$, such that
\mbox{$\val(\bp)=\oh(\psi(\bp))$}, $\bp\in \overline\bp^\infty$,
\item[(A4)] the sequence $\overline\bp$ is generic among the
sequences satisfying the above conditions.
\end{enumerate}}

Define the multiplicity function
$\mu:\overline\bp\cap(\K^*)^2\to\Z_{>0}$ such that:
\begin{itemize}\item for $\bp\in\overline\bp\cap(\K^*)^2$, $\val(\bp)=\ox\in
h(G^{\mr})$, put
\begin{equation}\mu(\bp)=\sum_{\gam\in G^{(\mr)},\
h(\gam)=\bx}\mt(\gam)\ ,\label{e12}\end{equation} \item for the
points $\bp_{1,\bx},\bp_{2,\bx}$, where
$\val(\bp_{1,\bx})=\val(\bp_{2,\bx})=\bx\in h(G^{(\dmp)})$, put
\begin{equation}\mu(\bp_{1,\bx})=m_1,\
\mu(\bp_{2,\bx})=m_2,\quad(m_1,m_2)= \sum_{\gam\in G^{(\dmp)},\
h(\gam)=\bx}\mt(\gam)\ .\label{e13}\end{equation}
\end{itemize} From this definition and from the count of the Euler
characteristic of $\oGam$, we derive
\begin{equation}\sum_{\bp\in\overline\bp\cap(\K^*)^2}
\mu(\bp)+|\overline\bp^\infty|-|\Gam^0_\infty|=g-1\
.\label{e19}\end{equation}

Let $\Del'\subset\R^2$ be a convex lattice polygon such that there
is another lattice polygon (or segment, or point) $\Del''$
satisfying $\Del'+\Del''=\Del$. Then we have a well defined line
bundle ${\cal L}_{\Del'}$ on $\Tor_\K(\Del)$. Let
$\overline\bp'\subset\overline\bp$ and
$\mu':\overline\bp'\cap(\K^*)^2\to\Z_{>0}$ be such that
$\mu'(\bp)\le \mu(\bp)$ for all $\bp\in\overline\bp'$. Let
$(\bet^\sig)'\in\Z_{\ge0}^\infty$, $\sig\subset\partial\Del$, be
such that $(\bet^\sig)'\le\bet^\sig$ for all sides $\sig$ of $\Del$.
We say that the tuple
$(\Del',g',\overline\bp',\mu',\{(\bet^\sig)'\}_{\sig\subset\partial\Del})$,
where $g'\in\Z_{\ge0}$, is compatible, if \begin{itemize}\item
$\|(\bet^\sig)'\|_1=\langle c_1({\cal
L}_{\Del'}),\Tor_\K(\sig)\rangle$ for all sides $\sig$ of $\Del$,
\item $(\bet^\sig)'_i\ge|\{\bp\in\overline\bp'\cap\Tor_\K(\sig)\ :\ \psi(\bp)=\gam\in
e\in\Gam^1_\infty,\ w(e)=i\}|$ for all sides $\sig$ of $\Del$ and
all $i=1,2,...$, \item $g'\le|\Int(\Del'\cap\Z^2)|$, and
$$\sum_{\bp\in\overline\bp'\cap(\K^*)^2}\mu'(\bp)+|\overline\bp'\cap\overline\bp^\infty|
-\sum_{\sig\subset\partial\Del}\|(\bet^\sig)'\|_0=g'-1\ .$$
\end{itemize} In view of (\ref{e19}) and
$|\Gam^0_\infty|=\sum_{\sig\subset\partial\Del}\|\bet^\sig\|_0$, the
tuple
$(\Del,g,\overline\bp,\mu,\{\bet^\sig\}_{\sig\subset\partial\Del})$
is compatible.

For any compatible tuple
$(\Del',g',\overline\bp',\mu',\{(\bet^\sig)'\}_{\sig\subset\partial\Del})$,
introduce the set ${\cal
C}(\Del',g',\overline\bp',m',\{(\bet^\sig)'\}_{\sig\subset\partial\Del})$
of reduced irreducible curves $C\in|{\cal L}_{\Del'}|$ passing
through $\overline\bp'$ and such that
\begin{itemize}\item the points
$\bp\in\overline\bp'\cap\overline\bp^\infty$ are nonsingular for
$C$, and
$$(C\cdot\Tor_\K(\partial\Del))_\bp=w(e)\ ,$$ where
$\gam=\psi(\bp)\in G_\infty$, and $e\in\Gam^1_\infty$ merges to
$\gam$,
\item the local branches of $C$ centered at the points of $C\cap\Tor_\K(\partial\Del)$ are
smooth, and, for each side $\sig$ of $\Del$ and each $i=1,2,...$,
there are precisely $(\bet^\sig)'_i$ local branches $P$ of $C$
centered at $C\cap\Tor_\K(\sig)$ such that
$$(P\cdot\Tor_\K(\sig))=i\ ,$$ $i=1,2,...$,
\item $C$ has genus $\le g'$,
\item at each point $\bp\in\overline\bp'\cap(\K^*)^2$, the multiplicity of $C$ is $\mt(C,\bp)\ge\mu'(\bp)$.
\end{itemize}

We now impose new conditions on the algebraic pathchworking data:

{\it\begin{enumerate}\item[(A5)] For any compatible tuple
$(\Del',g',\overline\bp',m',\{(\bet^\sig)'\}_{\sig\subset\partial\Del})$,
the set ${\cal
C}(\Del',g',\overline\bp',m',\{(\bet^\sig)'\}_{\sig\subset\partial\Del})$
is finite, all the curves $C\in{\cal
C}(\Del',g',\overline\bp',m',\{(\bet^\sig)'\}_{\sig\subset\partial\Del})$
are immersed, have genus $g'$, and have multiplicity
$\mt(C,\bp)=m'(\bp)$ at each point
$\bp\in\overline\bp'\cap(\K^*)^2$; furthermore, \begin{equation}
H^1(C^\nu,{\cal J}_Z(C^\nu))=0\ ,\label{e35}\end{equation} where
$C^\nu$ is the normalization and ${\cal J}_Z(C^\nu)$ is the (twisted
with $C^\nu$) ideal sheaf of the zero-dimensional scheme $Z\subset
C^\nu$ which contains the lift of $\overline\bp$ and of the points
of tangency of $C$ and $\Tor_\K(\partial\Del)$ upon $C^\nu$, and
which has length \mbox{$(C\cdot\Tor_\K(\partial\Del))_\bp$} at the
lift of $\bp\in\overline\bp\cap\Tor_\K(\partial\Del)$, and the
length \mbox{$(C\cdot\Tor_\K(\partial\Del))_\bz-1$} at the lift of
each point $\bz\in
C\cap\Tor_\K(\partial\Del)\backslash\overline\bp$.
\end{enumerate}}

Here we verify the condition (A5) for the versions of the
patchworking theorem used in \cite{IKS,Sh08}.

\begin{lemma}\label{l4} Condition (A5) holds if \begin{itemize}\item
either $\mu(\bp)=1$ for all $\bp\in\overline\bp\cap(\K^*)^2$,
\item or the surface $\Sig=\Tor_\K(\Del)$ is one of $\PP^2$, $\PP^2_k$
with $1\le k\le 3$, $(\PP^1)^2$, the configuration
$\overline\bp^\infty$ is contained in one toric divisor $E$ of
$\Sig$, and $$|\{\bp\in\overline\bp\cap(\K^*)^2\ :\
\mu(\bp)>1\}|\le\begin{cases} 4,\quad &\Sig=\PP^2,\\ 5-E^2-k,\quad
&\Sig=\PP^2_k,\\ 3\quad &\Sig=(\PP^1)^2\ .\end{cases}$$
\end{itemize}

Moreover, all the curves $C$ in the considered sets are nonsingular
along $\Tor_\K(\partial\Del)$, are nodal outside $\overline\bp$, and
have ordinary singularity of order $m'(\bp)$ at each point
$\bp\in\overline\bp'\cap(\K^*)^2$.
\end{lemma}

{\bf Proof.} We prove the statement only for the original data
$(\Del,g,\overline\bp,m,\{\bet^\sig\}_{\sig\subset\partial\Del})$,
since the other compatible tuples can be treated in the same way.

Observe that, in the first case, each curve $C\in{\cal
C}(\Del,g,\overline\bp,m,\{\bet^\sig\}_{\sig\subset\partial\Del})$
satisfies
\begin{equation}\sum_{\renewcommand{\arraystretch}{0.6}
\begin{array}{cc}
\scriptstyle{\bp\in\overline\bp\cap(\K^*)^2}\\
\scriptstyle{\mu(\bp)>1}
\end{array}}\mu(\bp)<|C\cap\Tor(\partial\Del)|-1\
.\label{e33}\end{equation} In the second situation, except for
finitely many lines or conics (which, of course, satisfy (A5)), the
other curves obey (\ref{e33}) by Bezout theorem (just consider
intersections with suitable lines or conics - we leave this to the
reader as a simple exercise).

So, we proceed further under the condition (\ref{e33}). Let
$\overline\bp'=\{\bp\in\overline\bp\cap(\K^*)^2\ :\ \mu(\bp)>1\}$.
Consider the family ${\cal C}'$ of reduced irreducible curves
$C'\in|{\cal L}_\Del|$ of genus at most $g$, which have multiplicity
$\ge \mu(\bp)$ at each point $\bp\in\overline\bp'$, whose local
branches centered along $\Tor_\K(\partial\Del)$ are nonsingular, and
the number of such branches crossing the toric divisor
$\Tor_\K(\sig)\subset\Tor_\K(\Del)$ with multiplicity $i$ equals
$\bet^\sig_i$ for all sides $\sig$ of $\Del$ and all $i=1,2,...$

The classical deformation theory argument (see, for instance,
\cite{GK,GL}) Zariski tangent space to ${\cal C}'$ at $C\in{\cal
C}:={\cal
C}(\Del,g,\overline\bp,m,\{\bet^\sig\}_{\sig\subset\partial\Del})$
is naturally isomorphic to $H^0(C^\nu,{\cal J}_Z(C^\nu))$, where
$C^\nu$, ${\cal J}_Z(C^\nu)$ are defined in (A5). So, we have
$$\deg
Z=\sum_{\bp\in\overline\bp'}\mu(\bp)+C\cdot\Tor_\K(\Del)-|C\cap\Tor(\partial\Del)|=\sum_{\bp\in\overline\bp'}\mu(\bp)
-CK_\Sig-|C\cap\Tor(\partial\Del)|$$ \begin{equation}<-CK_\Sig-1\
.\label{e31}\end{equation} Hence (see \cite{GK,GL}) $H^1(C^\nu,{\cal
J}_Z(C^\nu))=0$, which yields
$$h^0(C^\nu,{\cal J}_Z(C^\nu))=C^2-2\del(C)-\deg Z-g(C)+1$$
$$=-CK_\Sig+2g(C)-2-\deg Z-g(C)+1$$
\begin{equation}=g(C)-1+|C\cap\Tor(\partial\Del)|-\sum_{\bp\in\overline\bp'}\mu(\bp)\stackrel{\text{(\ref{e19})}}{=}
|\overline\bp\backslash\overline\bp'|-(g-g(C))\
.\label{e32}\end{equation} Since
$\overline\bp\backslash\overline\bp'$ is a configuration of generic
points (partly on $\Tor_\K(\partial\Del)$), we derive that $g(C)=g$
and that ${\cal C}$ is finite.

For the rest of the required statement, we assume that a curve
$C\in{\cal C}$ is either not nodal outside $\overline\bp$, or has
singularities on $\Tor_\K(\partial\Del)$, or has at some point
$\bp\in\overline\bp\cap(\K^*)^2$ a singularity more complicated than
an ordinary point of order $\mu(\bp)$. Then (cf. the argument in the
proof of \cite[Proposition 2.4]{OS}) one can find a zero-dimensional
scheme $Z\subset Z'\subset C^\nu$ of degree $\deg Z'=\deg Z+1$ such
that the Zariski tangent space to ${\cal C}'$ at $C$ is contained in
$H^0(C^\nu,{\cal J}_{Z'}(C^\nu))$. However, then one derives from
(\ref{e31}) that $\deg Z'<-CK_\Sig$, and hence again
$H^1(C^\nu,{\cal J}_{Z'}(C^\nu))=0$, which in view of (\ref{e32})
will lead to
$$h^0(C^\nu,{\cal
J}_{Z'}(C^\nu))=|\overline\bp\backslash\overline\bp'|-1\ ,$$ what
finally implies the emptiness of ${\cal C}$. \proofend

\subsection{Algebraic curves over $\K$ and tropical curves}
If $C\in|{\cal L}_\Del|$ is a curve on the toric surface
$\Tor_\K(\Del)$, then the closure
\mbox{$\Cl(\val(C\cap(\K^*)^2)))\subset\R^2$} supports an EPT-curve
$T$ with Newton polygon $\Del$ (cf. section \ref{sec4}) which is
defined by a convex piece-wise linear function (\ref{e2}) coming
from a polynomial equation $F(z)=0$ of $C$ in $(\K^*)^2$:
\begin{equation}F(z)=\sum_{\omega\in\Del\cap\Z^2}A_\omega z^\omega,\quad
A_\omega\in\K,\  c_\omega=\val(A_\omega),\quad z\in(\K^*)^2\
.\label{e30}\end{equation} The EPT-curve obtained does not depend on
the choice of the defining polynomial of $C$ and will be denoted by
$\Trop(C)$.

Observe also that the polynomial (\ref{e30}) can be written
$$F(z)=\sum_{\omega\in\Del\cap\Z^2}\left(a_\omega+O(t^{>0})
\right)t^{\nu_T(\omega)}z^\omega$$ with the convex piece-wise linear
function $\nu:\Del\to\R$ as in section \ref{sec4} and the
coefficients $a_\omega\in\C$ non-vanishing at the vertices of the
subdivision $S_T$ of $\Del$.

\subsection{Patchworking theorems} {\bf The algebraically closed
version.}

\begin{theorem}\label{t1} Given the patchworking data, a PPT-curve $\hat Q$ and a configuration
$\overline\bp$, satisfying all the conditions of section \ref{sec2},
there exists a subset ${\cal C}(\hat Q)\subset{\cal
C}(\Del,g,\overline\bp,m,\{\bet^\sig\}_{\sig\subset\partial\Del})$
of $M(Q)$ curves $C$ such that $\Trop(C)=h_*Q$. Furthermore, for any
distinct (non-isomorphic) curves $\hat Q_1$ and $\hat Q_2$, the sets
${\cal C}(\hat Q_1)$ and ${\cal C}(\hat Q_2)$ are disjoint.
\end{theorem}

\begin{remark}
We would like to underline one useful consequence of Theorem
\ref{t1}: The PPT-curve $Q$ and the multiplicities of its marked
curves must satisfy the restrictions known for the respective
algebraic curves with multiple points.
\end{remark}

\medskip
\noindent {\bf The real version.} In addition to all the above
hypotheses, we assume the following: {\it
\begin{enumerate}\item[(R6)] the configuration $\overline\bp$ is
$\Conj$-invariant,
$\val(\Re\overline\bp)\cap\val(\Im\overline\bp)=\emptyset$, where
\mbox{$\Re\overline\bp:=\Fix(\Conj\big|_{\overline\bp})$} and
$\Im\overline\bp=\overline\bp\backslash(\Re\overline\bp)$,
\item[(R7)] the PPT-curve $Q$ possesses a real structure
$\conj:Q\righttoleftarrow$\:, $G=\Re G\cup\Im G$ such that

(i) the bijection $\psi$ from (A3) takes $G_\infty\cap\Re G$ into
$\Re\overline\bp\cap\Tor_\K(\partial\Del)$ and takes
$G_\infty\cap\Im G$ into $\Im\overline\bp\cap\Tor_\K(\partial\Del)$,
respectively,

(ii) $h(G_o\cap\Re G)\subset\val(\Re\overline\bp)$, $h(G_0\cap\Im
G\cap\Gam^0)\subset\val(\Im\overline\bp)$,

(iii) $\Re G\cap G^{(\dmp)}\cap\Im\oGam=\emptyset$,

(iv) if $\gam\in G_0\cap\Im G\cap\Im\oGam$, $\mt(\gam)=(1,0)$, then
$\mt(\conj(\gam))=(0,1)$,

(v) if $e\in\Gam^1_\infty$, $e\subset\Re\oGam$, then $w(e)$ is odd.
\end{enumerate}}

\begin{theorem}\label{t2}
In the notations and hypotheses of Theorem \ref{t1} and under
assumptions (R1)-(R7), the following holds
\begin{equation}\sum_{C\in\Re{\cal C}(\hat Q)}W_{\Sig,\overline\bp}(C)=W(Q)\
,\label{e36}\end{equation} where $\Re{\cal C}(\hat Q)$ is the set of
the real curves in ${\cal C}(\hat Q)$, and
$W_{\Sig,\overline\bp}(C)$ defined in (\ref{e34}).
\end{theorem}

\subsection{Proof of Theorem \ref{t1}} Our argument is as follows.
First, we dissipate each multiple point $\bp\in\overline\bp$ of
multiplicity $k>1$ into $k$ generic simple points (in a neighborhood
of $\bp$), and then, using the known patchworking theorems
(\cite[Theorem 5]{Sh0} and \cite[Theorem 2.4]{Sh-g}\footnote{The
complete proof is provided in \cite{Sh-g}.}), obtain $M(Q)$ curves
$C\in|{\cal L}_\Del|$ of genus $g$ matching the deformed
configuration $\overline\bq$. After that, we specialize the
configuration $\overline\bq$ back into the original configuration
$\overline\bp$ and show that each of the constructed curves
converges to a curve with multiple points and tangencies as asserted
in Theorem \ref{t1}.

\begin{remark}
The deformation part of our argument works well in a rather more
general situation, whereas the degeneration part appears to be more
problematic, and at the moment we do not have a unified approach to
treat all possible degenerations which may lead to algebraic curves
with multiple points.
\end{remark}

Following \cite[Section 3]{Sh0}, we obtain the algebraic curves $C$
over $\K$ as germs of one-parameter families of complex curves
$C^{(t)}$, $t\in(\C,0)$, with irreducible fibres $C^{(t)}$, $t\ne
0$, of genus $g$ and a reducible central fibre $C^{(0)}$. The given
data of Theorem \ref{t1} provide us with a collection of suitable
central fibres $C^{(0)}$ out of which we restore the families using
the patchworking statement \cite[Theorem 5]{Sh0}.

\medskip

{\bf Step 1.} We start with a simple particular case which later
will serve as an element of the proof in the general situation.
Assume that $Q$ is a rational, simple, regularly end-marked
PPT-curve, $h_*Q\subset\R^2_\Del$ is a (compactified) nodal embedded
plane tropical curve, $\overline\bp\subset\Tor_\K(\partial\Del)$,
$G=G_\infty$, $\overline\bx=\overline\bx^\infty$, and
$\overline\bp\stackrel{\val}{\rightarrow}\overline\bx\stackrel{\oh}{\leftarrow}G$
are bijections. Here $M(Q)$ is given by formula (\ref{e17}), and
this number of required rational curves $C\subset\Tor_\K(\Del)$ is
obtained by a direct application of \cite[Theorem 5]{Sh0}.

The combinatorial part of the patchworking data for the construction
of curves over $\K$ consist of the tropical curve $Q$ which defines
a piece-wise linear function $\nu:\Del\to\R$ and a subdivision $S:\
\Del=\Del_1\cup...\cup\Del_N$ (see section \ref{sec4}). The
algebraic part of the patchworking data includes the limit curves
$C_k\subset\Tor(\Del_k)$, the deformation patterns $C_e$ associated
with the (finite length edges) edges $e\in\Gam^1$ (see \cite[Section
5.1]{Sh0} and \cite[Section 2.1]{Sh-g}), and the refined conditions
to pass through the fixed points (see \cite[Section 5.4]{Sh0} and
\cite[Section 2.5.9]{IMS}).

First, we orient the edges of $\Gam$ as in Lemma \ref{l2}(ii). Then
define complex polynomials $f_e$, $e\in\oGam^1$, and $f_v$,
$v\in\Gam^0$, in the following inductive procedure. In the very
beginning, for the $\Gam$-ends $e$ with marked points $\gam$, we
define
\begin{equation}f_e(x,y)=(\eta^q x^p-\xi^p y^q)^{w(e)}\
,\label{e14}\end{equation} where $\bu(e)=(p,q)$ and $(\xi,\eta)$ are
quasiprojective coordinates of the point $\ini(\bp)$ on
$\Tor(e)\subset\Tor(\Del)$ such that $\gam=\psi(\bp)$ (here $\psi$
is the bijection from condition (A3) above). Define a linear order
on $\Gam^0$ compatible with the orientation of $\Gam$. On each
stage, we take the next vertex $v\in\Gam^0$ and define $f_v$ and
$f_e$, where $e$ is the edge emanation from $v$. Namely, the
polynomials $f_{e_1},f_{e_2}$ associated with the two edges merging
to $v$, determine points $z_1\in\Tor(\sig_1)$, $z_2\in\Tor(\sig_2)$
on the surface $\Tor(\Del_v)$, $\sig_1,\sig_2$ being the sides of
$\Del_v$ orthogonal to $\oh(e_1)$, $\oh(e_2)$, respectively, and we
construct a polynomial $f_v$ with Newton polygon $\Del_v$ which
defines an irreducible rational curve $C_v\subset\Tor(\Del_v)$,
nonsingular along $\Tor(\partial\Del_v)$, crossing $\Tor(e_i)$ at
$z_i$, $i=1,2$, and crossing $\Tor(e)$ at one point $z_0$ (at which
one has $(C_v\cdot\Tor(e))_{z_0}=w(e)$). By \cite[Lemma 3.5]{Sh0},
up to a constant factor there are $|\Del_v|/(w(e_1)w(e_2))=M(Q,v)$
choices for such a polynomial $f_v$. After that we define $f_e(x,y)$
via (\ref{e14}) with $\xi,\eta$ the (quasihomogeneous) coordinates
of $z_0$ in $\Tor(\sig)$, where $\sig$ is the side of $\Del_v$
orthogonal to $\oh(e)$. So, the limit curves
$C_k\subset\Tor(\Del_k)$ are $C_v$ for the triangles $\Del_k$ dual
to $h(v)$, and are given by $f_{e_1}f_{e_2}$, where
$e_1,e_2\in\oGam^1$ appear in the decomposition (\ref{e3}) of a
parallelogram $\Del_k$.

The set of limit curves is completed by a set of deformation
patterns (see \cite[Sections 3.5 and 3.6]{Sh0}) as follows. Namely,
for each edge $e\in\Gam^1$ with $w(e)>1$, the deformation pattern is
an irreducible rational curve $C_e\subset\Tor(\Del_e)$, where
$\Del_e:=\conv\{(0,1),(0,-1),(w(e),0)\}$, whose defining (Laurent)
polynomial $f_e(x,y)$ has the zero coefficient of $x^{w(e)-1}$ and
the truncations to the edges $[(0,1),(w(e),0)]$ and
$[(0,-1),(w(e),0)]$ of $\Del_e$ fitting the polynomials $f_{v_1}$,
$f_{v_2}$, where $v_1,v_2$ are the endpoints of $e$ (see the details
in \cite[Sections 3.5 and 3.6]{Sh0}). Remind that, by \cite[Lemma
3.9]{Sh0} there are $w(e)=M(Q,e)$ suitable polynomials $f_e$.

The conditions to pass through a given configuration $\overline\bp$
do not admit a refinement. Indeed, following \cite[Section
5.4]{Sh0}, we can turn a given fixed point $\bp$ into $(\xi,0)$,
$\xi=\xi^0+O(t^{>0})\in\K$, by means of a suitable toric
transformation. Then, in \cite[Formula (6.4.26)]{Sh0}, the term with
the power $1/m$ will vanish.\footnote{The mentioned term contains
$\eta^0_s$, the initial coefficient of the second coordinate of
$\bp$, and not $\xi^0_s$ as appears in the published text. The
correction is clear, since in the preceding formula for $\tau$ one
has just $\eta^0_s$.}

The above collections of limit curves and deformation patterns
coincide with those considered in \cite{Sh0}, the transversality
hypotheses of \cite[Theorem 5]{Sh0} are verified in \cite[Section
5.4]{Sh0}. Hence, each of the
$$\prod_{v\in\Gam^0}M(Q,v)\cdot\prod_{e\in\Gam^1}M(Q,e)=M(Q)$$
above patchworking data gives rise to a rational curve
$C\subset\Tor_\K(\Del)$ as asserted in Theorem \ref{t1}. Notice
that all these curves are nodal by construction.

\medskip

{\bf Step 2}. Now we come back to the general situation and deform
the given configuration $\overline\bp$ into the following new
configuration $\overline\bq$.

Each point $\bp=(\xi t^a+...,\eta
t^b+...)\in\overline\bp\cap(\K^*)^2$ with multiplicity $\mu(\bp)>1$
(defined by (\ref{e12}) or (\ref{e13})) we replace by $\mu(\bp)$
generic points in $(\K^*)^2$ with the same valuation image
$\val(\bp)=(-a,-b)$ and the initial coefficients of the coordinates
close to $\xi,\eta\in\C^*$, respectively. Furthermore, we extend the
bijection $\psi$ from (A3) up to a map $\psi:\overline\bq\to G$ in
such a way that
\begin{itemize}\item $\val\big|_{\overline\bq}=\oh\circ\psi$, \item
each point $\gam\in G\backslash\Gam^0$ has a unique preimage, each
point $\gam\in G\cap\Gam^0$ has precisely two preimages,
\item if $\gam\in G^{(\dmp)}\backslash\Gam^0$, $h(\gam)=\bx$, then $\psi^{-1}(\gam)$ is
close to $\bp_{1,\bx}$ or to $\bp_{2,\bx}$ according as
$\mt(\gam)=(1,0)$ or $(0,1)$, \item if $\gam\in
G^{(\dmp)}\cap\Gam^0$, $h(\gam)=\bx$, then $\psi^{-1}(\gam)$
consists of two points, one close to $\bp_{1,\bx}$ and the other
close to $\bp_{2,\bx}$.\end{itemize}

Next we construct a set ${\cal C}'\subset{\cal
C}(\Del,g,\overline\bq,1,\{\bet^\sig\}_{\sig\subset\partial\Del})$
of $M(Q)$ curves with the tropicalization $h_*Q$. By Lemma \ref{l4},
they are irreducible, nodal, of genus $g$, and with specified
tangency conditions along $\Tor_\K(\partial\Del)$.

\medskip

{\bf Step 3.} Similarly to Step 1, we obtain the limit curves from
a collection of polynomials in $\C[x,y]$ associated with the edges
and vertices of the parameterizing graph $\Gam$ of $Q$:

(i) Let $\gam\in G\backslash\Gam^0$ lie on the edge $e\in\oGam^1$.
Then we associate with the edge $e$ a polynomial $f_e(x,y)$ given
by (\ref{e14}) with the parameters described in Step 1.

(ii) Let $\gam\in G_0$ be a (trivalent) vertex $v$ of $\Gam$. Then
$f_v(x,y)$ is a polynomial with Newton triangle $\Del_v$ (see
section \ref{sec4}) defining in $\Tor(\Del_v)$ a rational curve
$C_v\in|{\cal L}_{\Del_v}|$, which crosses each toric divisor of
$\Tor(\Del_v)$ at one point, where it is nonsingular, and which
passes through the two points $\ini(\psi^{-1}(\gam))$. Observe that
by \cite[Lemma 2.4]{Sh1}, up to a constant factor there are
precisely $|\Del_v|$ polynomials $f_v$ as above (though the
assertion and the proof of \cite[Lemma 2.4]{Sh1} are restricted to
the real case, it works well in the same manner in the complex case
regardless the parity of the side length of $\Del_v$).

(iii) Edges emanating from a vertex $v\in \Gam^0\cap G_0$ do not
contain any other point of $G$ due to the $\Del$-general position,
and we define polynomials $f_e$ for them by formula (\ref{e14})
where $\xi,\eta$ are the (quasihomogeneous) coordinates of the
intersection point of $C_v$ with $\Tor(\sig)$, $\sig$ being the side
of $\Del_v$ orthogonal to $\oh(e)$.

(iv) Pick a connected component $K$ of $\oGam\backslash G$ and
orient it as in Lemma \ref{l2}(ii). Then we inductively define
polynomials for the vertices and closed edges of $K$: In each stage
we define polynomials $f_v$ and $f_e$ for a vertex $v$ and a simple
closed edge $e$ emanating from $v$, whereas the polynomials $f_{e'}$
for all the edges $e'$ of $K$ merging to $v$ are given. Each of the
latter polynomials defines a point on $\Tor(\partial\Del_v)$, and
these points are distinct. We denote their set by $X$. Then we
choose a polynomial $f_v(x,y)$ with Newton triangle $\Del_v$
defining an irreducible rational curve $C_v\subset\Tor(\Del_v)$
which \begin{itemize}\item is nonsingular along
$\Tor(\partial\Del_v)$, \item crosses $\Tor(\partial\Del_v)$ at each
point $z\in X$ with multiplicity $w(e')$, where the edge
$e'\in\oGam^1$ merging to $v$ is associated with a polynomial
$f_{e'}$ which determines the point $z$, \item crosses
$\Tor(\partial\Del_v)\backslash X$ at precisely one point
$z_0$.\end{itemize} Notice that $z_0$ is the unique intersection
point of $C_v$ with the toric divisor
\mbox{$\Tor(\sig)\subset\Tor(\Del_v)$}, where $\sig$ is orthogonal
to $\oh(e)$, and $(C_v\cdot\Tor(\sig))_{z_0}=w(e)$. We claim that up
to a constant factor there are precisely $M(Q,v)$ polynomials $f_v$
as required.

The case of a trivalent vertex $v$ was considered in Step 1. In
general, observe that the set of the required curves is finite,
since we impose
$$(C_v\cdot\Tor(\partial\Del_v))-1=-C_vK_{\Tor(\Del_v)}-1$$
conditions on the rational curves $C_v\in|{\cal L}_{\Del_v}|$, and
the conditions are independent by Riemann-Roch. The cardinality of
this set does not depend neither on the choice of a generic
configuration of fixed points on $\Tor(\partial\Del_v)$, nor on
the choice of an algebraically closed ground field of
characteristic zero. Thus, we consider the field $\K$ and pick the
fixed points on $\Tor_\K(\partial\Del_v)$ so that the valuation
takes them injectively to a $\Del_v$-generic configuration in
$\partial\R^2_{\Del_v}$. Then the rule (M5) and the construction
in Step 1 provide $M(Q,v)$ curves as required. The fact that there
are no other curves in consideration follows from a slightly
modified Mikhalkin's correspondence theorem (for detail, see, for
instance, \cite{Sh08}).

We then define $f_e(x,y)$ via (\ref{e14}) with $\xi,\eta$ the
(quasihomogeneous) coordinates of $z_0$ in $\Tor(\sig)$.

Summarizing we deduce that the number of choices of the curves
$C_v$, $v\in\Gam^0$, and $C_e$, $e\in\oGam^1$, equals
$$\prod_{v\in\Gam^0}M(Q,v)\ .$$

\medskip

{\bf Step 4.} Now we define the limit curves, the deformation
patterns, and the refined conditions to pass through fixed points.

For each polygon $\Del_k$ of the subdivision $S$ of $\Del$, the
limit curve $C_k\subset\Tor(\Del_k)$ is defined by the product of
the constructed above polynomials $f_v,f_e$ corresponding to the
summands in the decomposition (\ref{e3}) of $\Del_k$.

The deformation pattern for each edge $e\in\Gam^1$ such that
$w(e)>1$ is defined in the way described in Step 1.

At last, the condition to pass through a given point
$\bq\in\overline\bq$ such that $\gam=\psi(\bq)\in G$ lies in the
interior of an edge $e\in\oGam^1$ with $w(e)>1$, admits a refinement
(see \cite[Section 5.4]{Sh0} and \cite[Section 2.5.9]{IMS}) which in
its turn is defined up to the choice of a $w(e)$-th root of unity,
where $e\in\oGam^1$ contains $\gam$.

So, the total number of choices we made up to now is
$$\prod_{v\in\Gam^0}M(Q,v)\cdot\prod_{e\in\Gam^1}M(Q,e)\cdot\prod_{\gam\in G}M(Q,\gam)=M(Q)\ .$$

\medskip

{\bf Step 5.} Let us verify the hypotheses of the patchworking
theorem from \cite{Sh0,Sh-g}.

First requirement to the limit curves (see \cite{Sh0}, conditions
(A), (B), (C) in section 5.1, or \cite{Sh-g}, conditions (C1), (C2)
in section 2.1) is ensured by the generic choice of $\ini(\bq)$,
$\bq\in\overline\bq$. Namely, the limit curves do not contain
multiple non-binomial components ({\it i.e.} defined by polynomials
with nondegenerate Newton polygons), any two distinct components of
any limit curve $C_k\subset\Tor(\Del_k)$ intersect transversally at
non-singular points which all lie in the big torus
$(\C^*)^2\subset\Tor(\Del_k)$, and, finally, the intersection points
of any component of a limit curve $C_k$ with $\Tor(\partial\Del_k)$
are non-singular.

The main requirement is the transversality condition  for the limit
curves and deformation patterns (see \cite[Section 5.2]{Sh0} and
\cite[Section 2.2]{Sh-g}), which is relative to the choice of an
orientation of the edges of the underlying tropical curve. In
\cite{Sh0,Sh-g}, one considers an orientation of edges of the
embedded plane tropical curve (cf. section \ref{sec4}), which in our
setting is just $h_*(Q)\subset\R^2$. Here we consider the
orientation of the edges of the connected components of
$\oGam\backslash G$ as defined in Lemma \ref{l2}(ii). Since this
orientation does not define oriented cycles and since the
intersection points of distinct components of any limit curve with
the toric divisors are distinct, the proof of \cite[Theorem 5]{Sh0}
and \cite[Theorem 2.4]{Sh-g} with the orientation of $\Gam$ is a
word-for-word copy of the proof with the orientation of $h(\Gam)$.
Moreover, comparing with \cite{Sh0,Sh-g}, here we impose extra
conditions to pass through the points $\ini(\bq)$,
$\bq\in\overline\bq$.

The deformation patterns are transversal in the sense of
\cite[Definition 5.2]{Sh0}, due to \cite[Lemma 5.5(ii)]{Sh0}, where
both the inequalities hold, since the deformation patterns are nodal
(\cite[Lemma 3.9]{Sh0}), and thus do not contribute to the left-hand
side of the inequalities, whereas their right-hand sides are
positive.

The transversality of a limit curve $C_k\subset\Sig_k:=\Tor(\Del_k)$
in the sense of \cite[Definition 5.1]{Sh0} means the triviality
({\it i.e.} zero-dimensionality) of the Zariski tangent space at
$C_k$ to the stratum in $|{\cal L}_{\Del_k}|$ formed by the curves
which split into the same number of rational components as $C_k$
({\it i.e.} the components of $C_k$ do not glue up when deforming
along such a stratum) each of them having the same number of
intersection points with $\Tor(\partial\Del_k)$ as the respective
component of $C_k$ and with the same intersection number, and such
that all but one of these intersection points are fixed. In other
words the conditions imposed on each of the components of $C_k$
determine a stratum with the one-point Zariski tangent space.
Indeed, the above fixation of intersection numbers of a component
$C'$ of $C_k$ and all but one intersection points of $C'$ with
$\Tor(\partial\Del_k)$ imposed $-C'K_{\Sig_k}-1$ conditions which
all are independent due to Riemann-Roch on $C'$.

Thus, \cite[Theorem 5]{Sh0} applies, and each of the $M(Q)$ refined
patchworking data constructed above produces a curve
$C\subset\Tor_\K(\Del)$ as asserted in Theorem \ref{t1}.

\medskip

{\bf Step 6.} Now we specialize the configuration $\overline\bq$
into $\overline\bp$ and prove that each of the curves $C\in{\cal
C}'$ constructed above tends (in an appropriate topology) to some
curve $\hat C\in|{\cal L}_\Del|$.

To obtain the required limits, we introduce a suitable topology.
Since the variation of $\overline\bq$ does not affect its valuation
image, the same holds for the (variable) curves $C'\in{\cal C}'$,
and hence one can fix once forever the function $\nu:\Del\to\R$.
Then, writing each coordinate of any point $\bq\in\overline\bq$ as
$X=t^{-\val(X)}\Psi_{\bq,X}(t)$ and each coefficient of the defining
polynomials of $C$ as $A_\omega=t^{\nu(\omega)}\Psi_\omega(t)$, and,
assuming (without loss of generality) that all the exponents of $t$
in the above coordinates and coefficients are integral, we deal with
the following topology in the space of the functions
$\psi_\omega(t)$ holomorphic in a neighborhood of zero: Take the
$C^0$ topology in each subspace consisting of the functions
convergent in $|t|\le\eps$ and then define the inductive limit
topology in the whole space.

So, we assume that the variation of $\overline\bq$ reduces to only
variation of $\ini(\bq)$, \mbox{$\bq\in\overline\bq\cap(\K^*)^2$},
whereas the reminders of the corresponding series in $t$ stay
unchanged.

To show that that the families of the curves $C\in{\cal C}'$ have
limits, we recall that their coefficients appear as solutions to a
system of analytic equations, which is soluble by the implicit
function theorem due to the transversality of the initial (refined)
patchworking data (cf. \cite{Sh0,Sh-g}). Thus, to confirm the
existence of the limits of the curves $C\in{\cal C}'$, it is
sufficient to show that the system of equations and the (refined)
patchworking data have limits and the latter limit is transverse. In
particular, we shall obtain that, in each coefficient
$A_\omega=t^{\nu(\omega)}\Psi_\omega(t)$, $\omega\in\Del\cap\Z^2$,
the factor $\Psi_\omega$ converges uniformly in the family.

We start with analyzing the specialization of limit curves. Since
the given tropical curve $Q$ stays the same, we go through the
curves $C_v$, $v\in\Gam^0$, and $C_e$, $e\in\oGam^1$. Clearly, the
curves $C_e$, $e\in\oGam^1$, keep their form (\ref{e14}) with the
parameters $\xi,\eta$ possibly changing as $\ini(\bq)$ tends to
$\ini(\bp)$, $\bp\in\overline\bp$. Similarly, the curves $C_v$
corresponding to the vertices $v\in\Gam^0$ of valency $3$ remain as
described in Step 1, {\it i.e.} nodal nonsingular along
$\Tor(\partial\Del_v)$, and crossing each toric divisor at one
point. Furthermore, the curves $C_v$ corresponding to the
non-special vertices $v\in\Gam^0$ of valency $>3$, remain as
described in Step 3, paragraph (iv), since the intersection points
of $C_v$ with the toric divisors which correspond to the edges of
$\oGam$, merging to $v$, do not collate and remain generic in the
specialization as they are not affected by possible collisions of
the points $\ini(\bq)$, $\bq\in\overline\bq$. So, let us consider
the case of a special vertex $v\in\Gam^0$. By (T4) $C_v$ cannot
split into proper components, and hence specializes into an
irreducible rational curve. Furthermore, the intersection points of
$C_v$ with the toric divisors which correspond to the special edges
may collate forming singular points, centers of several smooth
branches. So, finally, the transversality conditions for such a
curve reduce to the fact, that the Zariski tangent space at $C_v$ to
the stratum in $|{\cal L}_{\Del_v}|$ consisting of rational curves
with given intersection points along the two toric divisors which
are related to the oriented edges of $\oGam$ merging to $v$, is
zero-dimensional. This is precisely the same stratum conditions as
in Step 5, and the argument of Step 5 (Riemann-Roch on the rational
curve $C_v$) shows that all the $-C_vK_{\Tor(\Del_v)}-1$ conditions
defining the stratum in the Severi variety parameterizing the
rational curves in $|{\cal L}_{\Del_v}|$ are independent.

Next, we notice that by assumption (T4) the possible collision of
intersection points of $C_v$ with $\Tor(\partial\Del_v)$ concerns
only transverse intersection points ({\it i.e.} those which
correspond to edges of weight $1$), and hence does not affect
neither the deformation patterns, nor the refined conditions to pass
through $\overline\bp$. Thus, each of the curves $C\in{\cal C}'$
degenerates into some curve $\hat C\in|{\cal L}_{\Del}|$ which is
given by polynomial with coefficients
$A_\omega=t^{\nu(\omega)}\Psi_\omega(t)$, $\omega\in\Del$,
containing factors $\Psi_\omega$ convergent uniformly in some
neighborhood of $0$ in $\C$.

\begin{remark}\label{newr} (1) Observe that the genus of $\hat C$
does not exceed the genus of $C$.

(2) Notice also that there is no need to study refinements of
possible singular points appearing in the above collisions of the
intersection points of $C_v$ with $\Tor(\partial\Del_v)$. Indeed,
the number of the transverse conditions we found equals the number
of parameters - hence no any extra ramification is possible.
\end{remark}

\medskip

{\bf Step 7.} Next we show that, at each point $\bp\in\overline\bp$
with $\mu(\bp)>1$, the obtained curve $\hat C$ has $\mu(\bp)$ local
branches.

Considering the point $\bp\in(\K^*)^2$ as a family of points
$\bp^{(t)}\in(\C^*)^2$, $t\ne0$, we claim that the curves $\hat
C^{(t)}\subset\Tor(\Del)$ have $\mu(\bp)$ branches at $\bp^{(t)}$,
$t\ne 0$. For, we will describe how glue up the limit curves forming
$\hat C^{(0)}$ when $\hat C^{(0)}$ deforms into $\hat C^{(t)}$,
$t\ne0$. Our approach is to compare the above gluing with the gluing
of the limit curves in the deformation of $C^{(0)}$ into $C^{(t)}$,
$t\ne0$, where $C\in{\cal C}'$ passes through the configuration
$\overline\bq$, and the comparison heavily relies on the one-to-one
correspondence between the limit curves of $\hat C$ and $C\in{\cal
C}'$ established in Step 6.

\begin{figure}
\setlength{\unitlength}{1cm}
\begin{picture}(12,11)(0,0)
\thinlines\put(1,9){\line(0,-1){2}}\put(1,9){\line(-2,1){1}}
\put(1,9){\line(1,0){2}}\put(1,8.95){\line(1,0){2.1}}
\put(3,9){\line(1,1){1}}\put(3,9){\line(0,-1){1.5}}\put(3,7.5){\line(-1,-1){0.5}}\put(3,7.5){\line(1,0){1}}
\put(3.1,8.95){\line(1,1){1}}\put(3.1,8.95){\line(0,-1){1.5}}\put(3.1,7.45){\line(-1,-1){0.5}}\put(3.1,7.45){\line(1,0){1}}
\put(0.9,8.85){$\bullet$}\put(2.9,8){$\bullet$}\put(3.05,8){$\bullet$}
\put(0.6,8.8){$v$}\put(2,9.1){$e'$}\put(2,8.6){$e$}\put(2.6,8){$\gam'$}
\put(3.3,8){$\gam$}\put(2,6){\rm (a)}
\thinlines\put(7,9){\line(0,-1){2}}\put(7,9){\line(-2,1){1}}
\put(7,9){\line(1,0){3}}\put(7,8.95){\line(1,0){4}}\put(10,9){\line(0,-1){2}}\put(10,9){\line(1,1){1}}
\put(11,8.95){\line(0,-1){2}}\put(11,8.95){\line(1,1){1}}\put(6.95,8.85){$\bullet$}
\put(7.8,9.1){$e_1$}\put(7.8,8.6){$e_s$}\put(9,8.9){$\bullet$}\put(9,8.8){$\bullet$}
\put(9,9.2){$\gam_1$}\put(9,8.6){$\gam_s$}\put(6.6,8.8){$v$}\put(10.9,8.8){$\bullet$}
\put(9.9,8.9){$\bullet$}\put(9.7,9.2){$v_1$}\put(11.2,8.8){$v_s$}\put(8.5,6){\rm
(b)}
\thinlines\put(0.5,1){\vector(0,1){4}}\put(0.5,1){\vector(1,0){5}}
\put(1.5,1){\line(0,1){3}}\put(0.5,1){\line(1,1){1}}\put(0.5,1){\line(1,2){1}}
\put(1.5,2){\line(1,0){2}}\put(1.5,3){\line(1,0){1}}\dashline{0.2}(0.5,4)(1.5,4)
\thicklines\put(0.5,1){\line(1,3){1}}
\put(1.5,4){\line(1,-1){3}}\put(0.1,3.9){$s$}\put(4,0.6){$s+1$}\put(2,0){\rm
(c)} \thinlines\put(7,1){\vector(0,1){4}}\put(7,1){\vector(1,0){5}}
\put(8,2){\line(0,1){2}}\put(10,1){\line(0,1){1}}\put(7,1){\line(1,1){1}}\put(7,1){\line(1,2){1}}
\put(8,2){\line(1,0){2}}\put(8,3){\line(1,0){1}}\dashline{0.2}(7,4)(8,4)
\thicklines\put(7,1){\line(1,3){1}}
\put(8,4){\line(1,-1){3}}\put(6.6,3.9){$s$}\put(10.5,0.6){$s+1$}\put(8.5,0){\rm
(d)}
\end{picture}
\caption{Illustration to Step 7 of the proof of Theorem
\ref{t1}}\label{newf}
\end{figure}
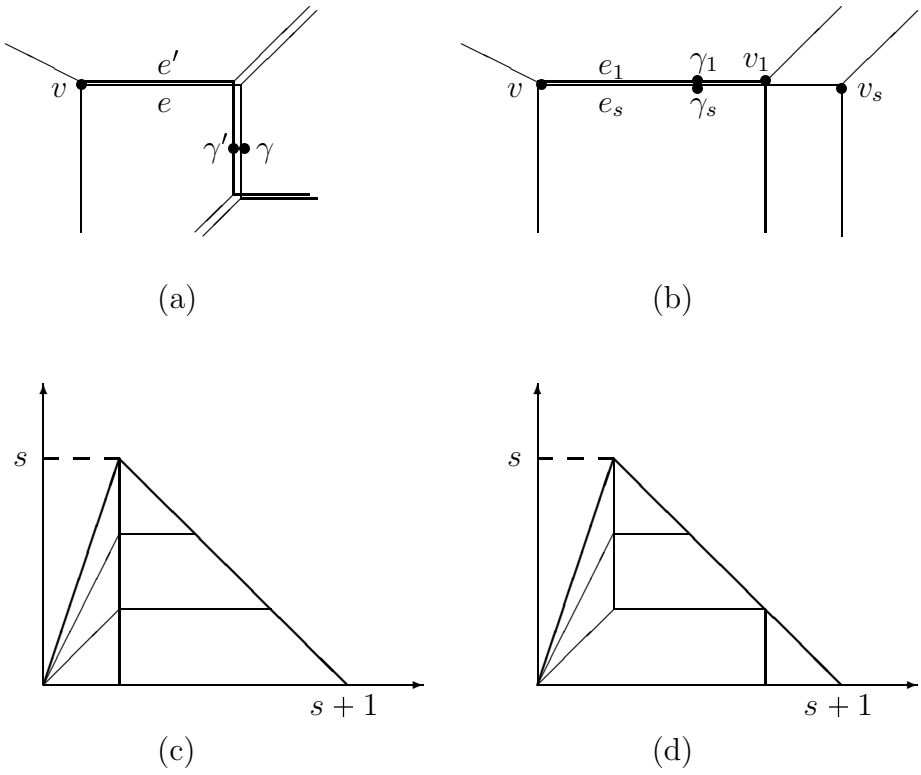

Let $\bq_1,...,\bq_s$ be all the points of the configuration
$\overline\bq$ which appear in the dissipation of the point
$\bp\in\overline\bp$ (cf. Step 2), and let $\gam_i=\psi(\bq_i)$,
$i=1,...,s$, be the corresponding marked points on $\Gam$ so that
$\gam_i\in e_i\in\Gam^1$, $i=1,...,s$. If the edges $h(e_i),h(e_j)$
intersect transversally at $V=h(\gam_i)=h(\gam_j)$, then $V$ is a
vertex of the plane tropical curve $h_*(Q)$ dual to a polygon
$\Del_V$ of the corresponding subdivision of $\Del$. The components
$C_i,C_j\subset\Tor(\Del_V)$ the curve $C^{(0)}$ passing through
$\ini(\bq_i),\ini(\bq_j)\in(\C^*)^2\subset\Tor(\Del_V)$,
respectively, intersect transversally in $(\C^*)^2$, and their
intersection points in $(\C^*)^2$ do not smooth up in the
deformation $C^{(t)}$, $t\ne 0$, and the same holds for the
corresponding components $\hat C_i,\hat C_j$ of $\hat C^{(0)}$
meeting at $\ini(\bp)\in(\C^*)^2\subset\Tor(\Del_V)$, since the
smoothing out of an intersection point $\ini(\bp)$ of $\hat C_i$ and
$\hat C_j$ would raise the genus of $\hat C$ above the genus of $C$
contrary to Remark \ref{newr}. Suppose that, in the above notation,
$h(e_i)$ and $h(e_j)$ lie on the same straight line, but $e_i,e_j$
have no vertex in common (see Figure \ref{newf}(a)). We consider the
case of finite length edges $e_i,e_j$; the case of ends can be
treated similarly. Let $v_i,v'_i$ be the vertices of $e_i$, and
$v_j,v'_j$ be the vertices of $e_j$. Their dual polygons
$\Del_{v_i},\Del_{v'_i},\Del_{v_j},\Del_{v'_j}$ (see section
\ref{sec4}) have sides $E_i,E'_i,E_j,E'_j$ orthogonal to $h(e_i)$.
In the deformation $C^{(0)}\to C^{(t)}$, $t\ne 0$, the limit curves
$C_i\subset\Tor(\Del_{v_i})$ and $C'_i\subset\Tor(\Del_{v'_i})$
passing through $\ini(\bq_i)\in\Tor(E_i)=\Tor(E'_i)$ glue up forming
a branch centered at $\bq_i^{(t)}$, and similarly the limit curves
$C_j\subset\Tor(\Del_{v_j})$ and $C'_j\subset\Tor(\Del_{v'_j})$
passing through $\ini(\bq_j)\in\Tor(E_j)=\Tor(E'_j)$ glue up forming
a branch centered at $\bq_j^{(t)}$. The same happens when $C$
specializes into $\hat C$, $\bq_i,\bq_j$ specialize into $\bp$,
since again the aforementioned restriction $g(\hat C)\le g(C)$ does
not allow the limit curves $\hat C_i,\hat C'_i$ glue up with the
limit curves $\hat C_j,\hat C'_j$.

The remaining case to study is given by the tropical data described
in condition (T7), section \ref{sec2}. Without loss of generality we
can assume that all the edges $e_1,...,e_s$ have a common vertex $v$
and their $h$-images lie on the same line (see an example in Figure
\ref{newf}(b)). Applying an appropriate invertible integral-affine
transformation, we can make the edges $e_1,...,e_s$ horizontal and
the point $\bx=\val(\bp)\in\R^2$ to be the origin. Correspondingly,
$v=(-\alp,0)$, $v_i=(\alp_i,0)$, $i=1,...,s$, with
$0<\alp_1\le...\le\alp_s\le\infty$ and
\begin{equation}\alp>\sum_{1\le i<s-1}\alp_i+2\alp_{s-1}\
.\label{newe7}\end{equation} In what follows we suppose that
$\alp_s<\infty$. The case $\alp_s=\infty$ admits the same treatment
as the case of finite $\alp_s\gg\alp$.

Let $\bq_i=\psi^{-1}(\gam_i)$, $1\le i\le s$, be the points of the
configuration $\overline\bq$ which appear in the deformation of the
point $\bp$ described in Step 2. Our assumptions yield that
$$\bp=(\xi+O(t^{>0}),\eta+O(t^{>0})),\quad\bq_i=(\xi_i+O(t^{>0}),\eta_i+O(t^{>0})),\
i=1,...,s\ ,$$ with some $\xi,\eta\in\C^*$, $\xi_i$ close to $\xi$,
$\eta_i$ close to $\eta$, $i=1,...,s$. Furthermore, the triangles
$\Del_v$ and $\Del_{v_i}$ dual to the vertices $v$ and $v_i$, $1\le
i\le s$, respectively, have vertical edges
$\sig\subset\partial\Del_v$ and $\sig_i\subset\partial\Del_{v_i}$,
$1\le i\le s$, along which the function $\nu$ (see section
\ref{sec4}) is constant. By assumptions (T4)-(T6), the limit curve
$C_v\subset\Tor(\Del_v)$ crosses the toric divisor $\Tor(\sig)$ at
the points $\eta_1,...,\eta_s$ with the total intersection
multiplicity $s$, and each of the limit curves
$C_{v_i}\subset\Tor(\Del_{v_i})$, $1\le i\le s$, crosses the toric
divisor $\Tor(\sig_i)$ at the unique point $\eta_i$ transversally,
and the corresponding limit curve $\hat C_{v_i}$ crosses
$\Tor(\sig_i)$ at the point $\eta$ transversally, too.

Now we move the points $\bq_1,...,\bq_s$ keeping their
$x$-coordinates and making \mbox{$(\bq_1)_y=...=(\bq_s)_y=(\bp)_y$}.
As shown in Step 6, the curve $C$ (depending on $\bq_1,...,\bq_s$)
converges to a curve $C'$ with the same Newton polygon, genus, and
tropicalization, and the limit curves of $C$ component-wise converge
to limit curves of $C'$. Consider now the polynomial $\widetilde
F(x,y):=F'(x,y+(\bp)_y)$, where the polynomial $F'(x,y)$ defines the
curve $C'$. As in the refinement procedure described in
\cite[Section 3.4]{Sh0} or \cite[Section 2.5.8]{IMS}, the
subdivision of the Newton polygon $\widetilde\Del$ of $\widetilde F$
contains the fragment bounded by the triangle
$\delta=\conv\{(0,0),(1,s),(s+1,0)\}$ (see Figure \ref{newf}(d))
which matches the points $\bq_1,...,\bq_s$. The corresponding
function $\widetilde\nu:\widetilde\Del\to\R$ takes the values
$$\widetilde\nu(0,0)=\alp,\quad\widetilde\nu(1,s)=0,\quad\widetilde\nu(k,s+1-k)=
\sum_{1\le i<k}\alp_i,\ k=2,...,s+1\ ,$$ along the incline part of
$\partial\delta$. The tropical limit of $\widetilde F$ restricted to
the above fragment consists of a subdivision of $\del$, determined
by some extension of the function $\widetilde\nu$ inside $\del$, and
of limit curves which must meet the following conditions:

- these limit curves glue up into a rational curve (with Newton
triangle $\del$), since, in the original tropical curve, the spoken
fragment corresponds to a tree (see Figure \ref{newf}(b));

- the intersection points $\bq$ of the curve $C_\delta:=\{\widetilde
F_\delta=0\}$ with the line $x=(\bp)_x$ such that $\val(\bq)_y\le0$
converge to $\bp$ as $\bq_1,...,\bq_s$ tend to $\bp$, where
$\widetilde F_\del$ is the sum of the monomials of $\widetilde F$
matching the set $\Del\cap\Z^2$, and the convergence is understood
in the topology of Step 6;

- the subdivision of $\del$ contains a segment $\widetilde\sig$ of
length $s$ lying inside the edge $[(0,0),(s+1,0)]$, along which the
function $\widetilde\nu$ is constant and such that the corresponding
toric divisor $\Tor(\widetilde\sig)$ intersects with the limit
curves at the points $\xi_1,...,\xi_s$.

These restrictions and inequality (\ref{newe7}) leave only one
possibility the subdivision of $\del$ shown in Figure
\ref{newf}(c,d) (the subdivision (c) for the case
$\alp>\alp_1+...+\alp_s$, and the subdivision (d) for the case
$\alp<\alp_1+...+\alp_s$). The limit curve
$C_{\del'}\subset\Tor(\del')$ for a triangle $\del'\subset\del$
having a horizontal base splits into $H(\del)$ distinct straight
lines (any of them crossing each toric divisor at one point), where
$H(\del')$ is the height. The limit curve
$\C_{\del'}\subset\Tor(\del')$ for a trapeze $\del'\subset\del$
splits into $H(\del')$ straight lines as above and the suitable
number of straight lines $x=\const$ (which reflect the splitting of
the trapeze into the Minkowski sum of a triangle with a horizontal
segment). All the limit curves are uniquely defined by the
intersections with the toric divisors $\Tor(\sig')$ for incline
segments $\sig'$ (in our construction, these data are determined by
the points $\overline\bq\backslash\{\bq_1,...,\bq_s\}$ and by the
condition $(\bq_1)_y=...=(\bq_s)_y=0$) and by intersections with
$\Tor(\widetilde\sig)$ introduced above. When $\bq_1,...,\bq_s$ tend
to $\bp$, the subdivision of $\del$ remains unchanged, whereas the
limit curves naturally converge component-wise. Then we immediately
derive that components of the limit curves passing through
$\ini(\bp)$ do not glue up together in the deformation $\hat
C^{(t)}$, $t\in(\C,0)$, since otherwise the (geometric) genus of
$\hat C^{(t)}$, $t\ne0$, would jump above the genus of $C$ which is
impossible (see Remark \ref{newr}).

\medskip

{\bf Step 8.} By assumption (A5), section \ref{sec2}, the curves
$\hat C\subset\Tor_\K(\Del)$ are immersed, irreducible, of genus
$g$, have multiplicity $\mu(\bp)$ at each point
$\bp\in\overline\bp\cap(\K^*)^2$, and satisfy the tangency
conditions with $\Tor_\K(\partial\Del)$ as specified in the
assertion of Theorem \ref{t1}. It remains to show that we have
constructed precisely $M(Q)$ curves $\hat C$.

Indeed, condition (\ref{e35}) implies that, for any dissipation of
each point $\bp\in\overline\bp\cap(\K^*)^2$ into $\mu(\bp)$ distinct
points there exists a unique deformation of $\hat C$ into a curve
$C\in{\cal C}$ such that a priori prescribed branches of $\hat C$ at
$\bp$ will pass through prescribed points of the dissipation.

Finally, we notice that the sets ${\cal C}(\hat Q_1)$ and ${\cal
C}(\hat Q_2)$ are disjoint for distinct (non-isomorphic PPT-curves
$\hat Q_1,\hat Q_2$. Indeed, the collections of limit curves as
constructed in Steps 1 and 3 appear to be distinct for distinct
curves $\hat Q_1$ and $\hat Q_2$ and the given configuration
$\overline\bp$. \proofend

\subsection{Proof of Theorem \ref{t2}} The curves $C\in\Re{\cal C}(\hat Q)$
constructed in the proof of Theorem \ref{t1} are immersed, and hence
the formula (\ref{e34}) for the Welschinger weight applies, thus the
left-hand side of (\ref{e36}) is well defined.

Next we go through the proof of Theorem \ref{t1} counting the
contribution to the right-hand side of (\ref{e36}).

First, we deform the configuration $\overline\bp$ as described in
Step 2, assuming that the deformed configuration $\overline\bq$ is
$\Conj$-invariant and that the map $\psi:\overline\bq\to G$ sends
$\Re\overline\bq=\overline\bq\cap\Fix(\Conj)$ to $\Re G$ and sends
$\Im\overline\bq=\overline\bq\backslash\Re\overline\bq$ to $\Im G$,
respectively. In particular, if $\bp\in\Re\overline\bp$, and the
points $\val(\bp)$ is an image of $r$ points of $\Re G\cap\Re\oGam$
and $s$ pairs of points of $\Re G\cap\Im\oGam$, then $\bp$ deforms
into $r$ real points and $s$ pairs of imaginary conjugate points.

Notice that the replacement of $\overline\bp$ by $\overline\bq$
causes a change of sign in the left-hand side of (\ref{e36}) and of
the quantity of the real curves in count in the right-hand side of
(\ref{e36}). Right now we explain the change of sign: the
dissipation of a real point $\bp$ as in the preceding paragraph
means that, for each curve $C\in{\cal C}'$ we count $r$ real
solitary nodes more in a neighborhood of $\bp$, since in the
non-deformed situation, the point $\bp$ should be blown up for the
computation of the Welschinger sign. This change is reflected in the
sign $(-1)^{\ell_1}$, $\ell_1=|\Re G\cap\Im\oGam|/2$, in the
right-hand side of formula (\ref{e16}).

Next we follow the procedure in Steps 3 and 4 of the proof of
Theorem \ref{t1} and construct $\Conj$-invariant collections of
limit curves, deformation patterns, and refined conditions to pass
through fixed points: \begin{itemize}\item by \cite[Proposition
8.1(i)]{Sh0} the existence of an even weight edge $e\in\Gam^1$,
$e\subset\Re\oGam$ annihilates the contribution to the Welschinger
number, and hence by (R7)(v) we can assume that all the edges
$e\subset\Re\oGam$ have odd weight, in particular, with the finite
length edges $e\subset\Re\oGam$ one can associate a unique real
deformation pattern with an even number of solitary nodes (cf.
\cite[Lemma 2.3]{Sh1}),
\item the limit curves associated with the vertices of $\Re\oGam$
contribute as designated in rules (W2)-(W4) in section \ref{sec1}
(cf. \cite[Lemmas 2.3, 2.4, and 2.5]{Sh1}),
\item the construction of limit curves and deformation patterns
associated with the vertices and edges of $\Im\oGam'$ (a half of
$\Im\oGam$), contributes as designated in rules (W1)\;-(W3) (cf.
with the complex formulas in the proof of Theorem \ref{t1} and with
\cite[Section 2.5]{Sh1}), accordingly, the data associated with
$\Im\oGam''$ are obtained by the conjugation,
\item the refinement of the condition to pass through fixed points
contributes as designated in rule (W2) as we have a unique
refinement for $\gam\in\Re\oGam$ and $w(e)$ refinements for $\gam\in
e\in\Im\oGam^1$.
\end{itemize} The remaining step is to explain the factor $2^{\ell_2}$,
$\ell_2=(|\Im G\cap\Im\oGam|-b_0(\Im\oGam))/2$, in formula
(\ref{e16}). Indeed, when constructing the limit curves associated
with the vertices of $\Im\oGam'$, we start with the respective fixed
points which all are imaginary in the configuration $\overline\bq$,
and thus we choose a point in each of the
$|G\cap\Im\oGam'|=|G\cap\Im\oGam|/2$ pairs of the corresponding
points in $\overline\bq$. Observe that, in the degeneration
$\overline\bq\to\overline\bp$, $|\Re G\cap\Im\oGam'|$ pairs of
imaginary points of $\overline\bq$ merge to real points in
$\overline\bp$, which leaves only $|\Im G\cap\Im\oGam|/2$ choices in
the original configuration $\overline\bp$. After all, we factorize
by the interchange of the components of $\Im\oGam$, coming to the
required factor $2^{\ell_2}$. \proofend

{\ncsc School of Mathematical Sciences \\[-21pt]

Raymond and Beverly Sackler Faculty of Exact Sciences\\[-21pt]

Tel Aviv University \\[-21pt]

Ramat Aviv, 69978 Tel Aviv, Israel} \\[-21pt]

{\it E-mail address}: {\ntt shustin@post.tau.ac.il}


\begin{thebibliography}{99}

\bibitem{GM1} Gathmann, A., and Markwig, H.: The numbers of tropical plane
curves through points in general position. {\it J. reine angew.
Math.} {\bf 602} (2007), 155--177.

\bibitem{GK} Greuel, G.-M., and Karras, U.: {\it
Families of varieties with prescribed singularities}. Compos. math.
{\bf 69} (1989), no. 1, 83--110.

\bibitem{GL} Greuel, G.-M., and Lossen, C.: {\it Equianalytic and equisingular families of
curves on surfaces}. Manuscripta math. {\bf 91} (1996), no. 3,
323--342.

\bibitem{IKS2} Itenberg, I. V.,
Kharlamov, V. M., and Shustin, E. I.: Logarithmic equivalence of
Welschinger and Gromov-Witten invariants. {\it Russian Math.
Surveys} {\bf 59} (2004), no. 6, 1093--1116.

\bibitem{IKS3} Itenberg, I. V.,
Kharlamov, V. M., and Shustin, E. I.: New cases of logarithmic
equivalence of Welschinger and Gromov-Witten invariants. {\it
Proc. Steklov Math. Inst.} {\bf 258} (2007), 65-73.

\bibitem{IKS} Itenberg, I., Kharlamov, V., and Shustin, E.: {\it
Recursive formulas and logarithmic asymptotics of Welschinger
invariants of real non-toric Del Pezzo surfaces}, in preparation.

\bibitem{IMS} Itenberg, I., Mikhalkin, G., and Shustin, E.: {\it Tropical algebraic
geometry}/ Oberwolfach seminars, vol. 35. Birkhauser, 2007.

\bibitem{Mi1} Mikhalkin, G.: Decomposition into pairs-of-pants for complex
algebraic hypersurfaces. {\it Topology} {\bf 43} (2004),
1035--1065.

\bibitem{Mi} Mikhalkin,~G.: Enumerative tropical algebraic
geometry in $\R^2$. {\it J. Amer. Math. Soc.} {\bf 18} (2005),
313--377.

\bibitem{NS} Nishinou, T., and Siebert, B.: Toric degenerations of toric
varieties and tropical curves. {\it Duke Math. J.} {\bf 135}
(2006), no. 1, 1--51.

\bibitem{OS} Orevkov, S., and Shustin, E.: {\it Pseudoholomorphic, algebraically
unrealizable curves}. Moscow Math. J. {\bf 3} (2003), no. 3,
1053--1083.

\bibitem{RST} Richter-Gebert, J., Sturmfels, B., and
Theobald, T.: First steps in tropical geometry. {\it Idempotent
mathematics and mathematical physics}, Contemp. Math. {\bf 377},
Amer. Math. Soc., Providence, RI, 2005, pp. 289--317.

\bibitem{Sh0} Shustin,~E.: A tropical approach to enumerative geometry. {\it Algebra i
Analiz} {\bf 17} (2005), no. 2, 170--214 (English translation: St.
{\it St. Petersburg Math. J.} {\bf 17} (2006), 343--375).

\bibitem{Sh91} Shustin,~E.: On manifolds of singular algebraic curves. {\it Selecta Math.
Sov.} {\bf 10}, no. 1, 27--37 (1991).

\bibitem{Sh1} Shustin,~E.: A tropical calculation of the Welschinger invariants
of real toric Del Pezzo surfaces. {\it J. Alg. Geom.} {\bf 15}
(2006), no. 2, 285--322 (corrected version at arXiv:math/0406099).

\bibitem{Sh-g} Shustin,~E.: Patchworking construction in the tropical enumerative
geometry. {\it Singularities and Computer Algebra}, C. Lossen and G.
Pfister, eds., Lond. Math. Soc. Lec. Notes Ser. {\bf 324}, Proc.
Conf. dedicated to the 60th birthday of G.-M. Greuel, Cambridge
Univ. Press, 2006, pp. 273--300.

\bibitem{Sh2} Shustin,~E.: Welschinger invariants of toric Del Pezzo
surfaces with non-standard real structures. {\it Proc. Steklov
Math. Inst.} {\bf 258} (2007), 219--247.

\bibitem{Sh08} Shustin,~E.: {\it New enumerative invariants and
correspondence theorems for plane tropical curves}, in preparation.

\bibitem{Vi2} Viro, O. Ya.:
Gluing of plane real algebraic curves and construction of curves
of degrees $6$ and $7$. {\it Lect. Notes Math.} {\bf 1060},
Springer, Berlin etc., 1984, pp. 187--200.

\bibitem{Vi3} Viro, O. Ya.:
Real algebraic plane curves: constructions with controlled
topology. {\it Leningrad Math. J.} {\bf 1} (1990), 1059--1134.

\bibitem{Vi4} Viro, O. Ya.: {\it Patchworking Real Algebraic Varieties}.
Preprint at arXiv:math/0611382.

\bibitem{Vi} Viro,~O.: Dequantization of Real Algebraic Geometry on a Logarithmic Paper.
{\it Proceedings of the 3rd European Congress of Mathematicians},
Birkh\"{a}user, Progress in Math. {\bf 201}, (2001), 135--146.

\bibitem{W1} Welschinger,~J.-Y.:
Invariants of real symplectic 4-manifolds and lower bounds in real
enumerative geometry. {\it Invent. Math.} {\bf 162} (2005), no. 1,
195--234.

\end{thebibliography}
\end{document}